\documentclass[preprint,12pt]{elsarticle}

\usepackage{amssymb}

\usepackage[utf8]{inputenc}
\usepackage{lmodern} 
\usepackage[T1]{fontenc}
\usepackage{microtype}

\usepackage{graphicx}
\usepackage{float}
\usepackage[hypcap]{caption}
\usepackage{subcaption}
\usepackage{algorithm2e}
\usepackage{array}
\usepackage{algorithmic}

\usepackage{url}

\usepackage{textcomp}
\usepackage{booktabs}
\usepackage{array}

\usepackage{amsfonts}

\usepackage{geometry}
\usepackage{tikz}
\usetikzlibrary{trees,positioning,shapes}

\usepackage{mathtools}
\usepackage{amssymb}
\usepackage{amsthm}
\usepackage{amsmath}
\usepackage{xcolor, colortbl}
\usepackage{appendix}

\usepackage{soul}
\usepackage{pgfplots,tikz}
\usepackage[eulergreek]{sansmath}
\usepackage{listings}
\usepackage{minted}
\pgfplotsset{compat=1.18} 
\usepackage[english]{babel}
\newcommand{\Prod}[2]{(#1, #2)}

\protect

\protect

\usepackage[colorlinks=true, linkcolor=blue, urlcolor=blue]{hyperref}

\journal{Computers and Mathematics with Applications}

\begin{document}

\begin{frontmatter}



\title{CO$_2$ sequestration hybrid solver using isogeometric alternating-directions and collocation-based robust variational physics informed neural networks (IGA-ADS-CRVPINN)}


\author[inst1]{Askold Vilkha}
\author[inst1]{Tomasz S\l{}u\.zalec}
\author[inst1]{Marcin \L{}o\'s}
\author[inst1]{Maciej Paszy\'nski\corref{cor1}}
\cortext[cor1]{Corresponding author. Email address: paszynsk@agh.edu.pl}

\affiliation[inst1]{organization={AGH University of Krakow, Faculty of Computer Science},%Department and Organization
            addressline={Al. Mickiewicza 30}, 
            city={Krak\'ow},
            postcode={30-059}, 
            country={Poland}}

\begin{abstract}
    This paper presents the hybrid solver for a $CO_2$ sequestration problem. The solver uses the IGA-ADS (IsoGeometric Analysis Alternating Directions solver) to compute the saturation scalar field update using the explicit method, and CRVPINN (Collocation-based Robust Variational Physics Informed Neural Networks solver) to compute the pressure scalar field. The study focuses on simulating the physical behavior of $CO_2$ in porous structures, excluding chemical reactions. The mathematical model is based on Darcy's Law. 
    The CRVPINN is pretrained on the initial pressure configuration, and the time step pressure updates require only 100 iterations of the Adam method per time step. We compare our hybrid IGA-ADS solver, coupled with the CRVPINN method, with a baseline of the IGA-ADS solver coupled with the MUMPS direct solver. Our hybrid solver is over 3 times faster on a single computational node from the ARES cluster of ACK CYFRONET.
 Future work includes extensive testing, inverse problem solving, and potential application to $H_2$ storage problems.
 \end{abstract}

\begin{keyword}
$CO_2$ sequestration \sep Isogeometric finite element method \sep Alternating-directions sovler \sep Physics Informed Neural Networks \sep Robust loss \sep Collocation method
\end{keyword}

\end{frontmatter}

\section{Introduction}
Given the significant influence of carbon dioxide ($CO_2$) on global warming, the development of carbon capture strategies is of the highest importance.
One of the prevalent strategies is the process called $CO_2$ sequestration.
The main idea of this process is to directly capture the $CO_2$ from the emission sources, such as power plants and factories, and store it in underground porous structures of both natural and human-made origin (abandoned oil wells, mines, aquifers) \cite{lackner_guide_2003,holloway_storage_2001,pruess_multiphase_2002}.
This strategy is more beneficial than others (ocean storage, mineralization, direct air capture), both economically and ecologically.
It is expected to be the most useful in the short and medium-term and can help mitigate the economic effects until more environmentally-friendly alternatives for energy production come into wide usage \cite{massarweh_co2_2024,bielinski_numerical_2007}. In the context of the $CO_2$ sequestration problem, four common situations are outlined: structural, capillary, dissolution, and mineralization \cite{massarweh_co2_2024}. 
We develop a solver that can be applied to the first three of them, as it does not take into account the chemical reactions in the reservoir, but can effectively simulate the physical behavior of the pumped gas in the porous structure, and in the case of pumping into the liquid brine solution.

In this work, we are implementing the $CO_2$ sequestration problem in the IGA-ADS (IsoGeometric Analysis using Alternating Directions Solver) software \cite{los_iga-ads_2017}.
This framework is a Finite-Element Method solver employing a sequence of isogeometric $L^2$ projections for solving time-dependent problems with higher-order and continuity B-spline-based discretization. It has been successfully applied to problems like heat transfer \cite{los_iga-ads_2017}, Maxwell's equations solver \cite{los_fast_2023}, and tumor growth simulations \cite{noauthor_tuning_2018}.
In this paper, the IGA-ADS solver is applied to compute the saturation scalar field.
Unfortunately, we also need to compute the pressure scalar field, which is not compatible with the IGA-ADS method, as it is described by a stationary PDE. Thus, we couple the IGA-ADS solver with the Collocation-Based Robust Variational Physics Informed Neural Networks (CRVPINN) \cite{los_collocation-based_2025}. The CRVPINN is an extension of the Robust Variational PINN (RVPINN) \cite{ROJAS2024116904}, 
introducing the collocation points and the Kronecker delta test functions to avoid expensive integrations of the weak residuals of RVPINN. The RVPINN uses the robust loss function. It modifies the weak residuals of the Variational PINN (VPINN)  \cite{KHARAZMI2021113547} by incorporating the inverse of the Gram matrix computed with the inner product, for which the weak formulation is bounded and inf-sup stable. The VPINN replaces the strong residuals of the Physics Informed Neural Network (PINN)  \cite{karniadakis_physics-informed_2021} with the weak residuals. In this sense, the CRVPINN method can be seen as a PINN method extended to discrete weak residuals, computed with the Kronecker delta test function, modified by the inverse of the Gram matrix, which makes the discrete weak residual-based loss function robust.

The hybrid IGA-ADS solver coupled with the CRVPINN solver is over 3 times faster than the IGA-ADS solver coupled with the MUMPS direct solver \cite{MUMPS:1,MUMPS:2}.
We believe that our method will provide a fast and reliable solver for this problem.
Moreover, in our future work, we also plan to experiment with the (hydrogen) $H_2$ storage problems and implement the inverse problem solver.

\section{CO2 sequestration problem formulation}
Following the previous works \cite{pruess_multiphase_2002,ebigbo_co2_2007,shokouhi_physics-informed_2021}, we use the analytical expression for the $CO_2$ sequestration problem (Darcy's Law).
\begin{equation}{\label{eq:1}}
    \phi \partial_t S_w = \nabla \cdot \left( \frac{S_w}{\mu_w} K \left( \nabla p - \rho_w \vec{g} \right) \right) + q_w
\end{equation}
\begin{equation}{\label{eq:2}}
    \phi \partial_t S_g = \nabla \cdot \left( \frac{S_g}{\mu_g} K \left( \nabla p - \rho_g \vec{g} \right) \right) + q_g
\end{equation}
where $S_w$ and $S_g$ are the brine solution and gas saturations, respectively, $\phi$ is the porosity, $\mu_w$ and $\mu_g$ are the viscosities of the brine solution and gas, $K$ is the permeability, $p$ is the pressure, $\rho_w$ and $\rho_g$ are the densities of the brine solution and gas, $g$ is the gravitational acceleration, and $q_w$ and $q_g$ are the source terms for the brine solution and gas, respectively.
A trivial constraint that $S_w + S_g = 1$ is also applied, so we will only compute one of the fields.
Taking this into account, we can rewrite the equation \eqref{eq:1} as follows:
\begin{equation}{\label{eq:3}}
    -\phi \partial_t S_g = \nabla \cdot \left( \frac{1 - S_g}{\mu_w} K \left( \nabla p - \rho_w \vec{g} \right) \right) + q_w
\end{equation}
Hence, we have a system of two equations \eqref{eq:3} and \eqref{eq:2} with two unknowns $S_g$ and $p$. 
To simplify the system even further, we can sum these equations and get the following equation:
\begin{equation}{\label{eq:4}}
    \nabla \cdot \left( \left( \frac{1 - S_g}{\mu_w} + \frac{S_g}{\mu_g} \right) K \nabla p \right) = \nabla \cdot \left( \vec{g} K\left( \frac{1 - S_g}{\mu_w} \rho_w + \frac{S_g}{\mu_g} \rho_g \right)\right) - (q_w + q_g)
\end{equation}
We will then use this equation to compute pressure $p$ after setting some initial saturation $S_g^0$. In each following time step, we compute pressure $p^n$ based on the saturation $S_g^n$ as follows:
\begin{equation}{\label{eq:5}}
    \nabla \cdot \left( \left( \frac{1 - S_g^n}{\mu_w} + \frac{S_g^n}{\mu_g} \right) K \nabla p^n \right) = \nabla \cdot \left(\vec{g} K \left( \frac{1 - S_g^n}{\mu_w} \rho_w + \frac{S_g^n}{\mu_g} \rho_g \right)\right) - (q_w^n + q_g^n)
\end{equation}
After computing $p^n$ and $S_g^n$, we can compute $S_g^{n+1}$ based on the pressure $p^n$ as follows (using the forward Euler method):
\begin{equation*}
    \phi \frac{S_g^{n+1} - S_g^n}{\tau} = \nabla \cdot \left( \frac{S_g^n}{\mu_g} K \left( \nabla p^n - \rho_g \vec{g} \right) \right) + q_g^n
\end{equation*}
\begin{equation}{\label{eq:6}}
    S_g^{n+1} = S_g^n + \phi^{-1} \tau \nabla \cdot \left( \frac{S_g^n}{\mu_g} K \left( \nabla p^n - \rho_g \vec{g} \right) \right) + \phi^{-1} \tau q_g^n
\end{equation}

%The programming implementation of the solver is done in the C++ programming language and can be found at the following link: \href{https://github.com/askoldvilkha/iga-ads-vilkha/tree/881587b8021a7a05c339ce9fde90b3d6acb24f62/examples/co2_sequestration}{GitHub}
%\footnote{https://github.com/askoldvilkha/iga-ads-vilkha/tree/develop/examples/co2\_sequestration}.

\section{Weak formulation of the problem}

First, let us consider the weak formulation of the saturation equation \eqref{eq:6}.
The computationally challenging part of this equation is the divergence term. 
Hence, we will derive the weak formulation to get rid of this operator. In order to do this, we multiply the equation by a test function $v$ and integrate over the computational domain $\Omega$:
\begin{equation}{\label{eq:7}}
    \int_{\Omega} S_g^{n+1} v \, d\Omega = 
    \int_{\Omega} (S_g^n  +  \phi^{-1} \tau q_g^n) v \, d\Omega + 
    \phi^{-1} \tau \int_{\Omega} \nabla \cdot \left( \frac{S_g^n}{\mu_g} K \left( \nabla p^n - \rho_g \vec{g} \right) \right) v \, d\Omega 
\end{equation}

Now, we can solve the last term using integration by parts:
\begin{equation*}
\begin{aligned}
    &\int_{\Omega} \nabla \cdot \left( \frac{S_g^n}{\mu_g} K \left( \nabla p^n - \rho_g \vec{g} \right) \right) v \, d\Omega = \\
    &\int_{\Gamma} \left( \frac{S_g^n}{\mu_g} K \left( \nabla p^n - \rho_g \vec{g} \right) \cdot \vec{n} \right) v \, d\Gamma - 
    \int_{\Omega} \left( \frac{S_g^n}{\mu_g} K \left( \nabla p^n - \rho_g \vec{g} \right) \right) \cdot \nabla v \, d\Omega
\end{aligned}
\end{equation*}
, where $\Gamma$ is the boundary of the computational domain $\Omega$, and $\vec{n}$ is the outward normal vector to the boundary $\Gamma$.
By selecting boundary conditions such that the boundary integral vanishes (in our case, we use zero Dirichlet boundary conditions on saturation $S_g$), the last term simplifies to:
\begin{equation}{\label{eq:8}}
    \int_{\Omega} \nabla \cdot \left( \frac{S_g^n}{\mu_g} K \left( \nabla p^n - \rho_g \vec{g} \right) \right) v \, d\Omega = 
    - \int_{\Omega} \left( \frac{S_g^n}{\mu_g} K \left( \nabla p^n - \rho_g \vec{g} \right) \right) \cdot \nabla v \, d\Omega
\end{equation}
Now, considering that vector $\vec{g}$ is pointed downwards (in the vertical $y$-axis), we can rewrite the last term as follows:
\begin{equation}{\label{eq:9}}
    \int_{\Omega} \nabla \cdot \left( \frac{S_g^n}{\mu_g} K \left( \nabla p^n - \rho_g \vec{g} \right) \right) v \, d\Omega =
    \int_{\Omega} \left( \frac{S_g^n}{\mu_g} K (\rho_g g v_y - \nabla p^n \cdot \nabla v) \right) \, d\Omega
\end{equation}
Hence, our weak formulation of the saturation equation \eqref{eq:6} becomes:
\begin{equation}{\label{eq:10}}
    \int_{\Omega} S_g^{n+1} v \, d\Omega = 
    \int_{\Omega} S_g^n v \, d\Omega + 
    \phi^{-1} \tau \int_{\Omega} \left( \frac{S_g^n}{\mu_g} K (\rho_g g v_y - \nabla p^n \cdot \nabla v) + q_g^n v\right) \, d\Omega
\end{equation}
Next, we will derive the weak formulation of the pressure equation \eqref{eq:5}.
Similar to the previous case, we multiply the equation by a test function $v$ and integrate over the computational domain $\Omega$:
\begin{equation}{\label{eq:11}}
\begin{aligned}
    &\int_{\Omega} \nabla \cdot \left( \left( \frac{1 - S_g^n}{\mu_w} + \frac{S_g^n}{\mu_g} \right) K \nabla p^n \right) v \, d\Omega = \\
    &\int_{\Omega} \nabla \cdot \left( \vec{g} K \left( \frac{1 - S_g^n}{\mu_w} \rho_w + \frac{S_g^n}{\mu_g} \rho_g \right) \right) v \, d\Omega - 
    \int_{\Omega} (q_w^n + q_g^n) v \, d\Omega
\end{aligned}
\end{equation}
As we have already established in the previous case, we can transform the divergence terms using the integration by parts similar to equation \eqref{eq:8}.
Thus, we get the following weak formulation of the pressure equation \eqref{eq:5}:
\begin{equation}{\label{eq:12}}
\begin{aligned}
    &\int_{\Omega} \left( \left( \frac{1 - S_g^n}{\mu_w} + \frac{S_g^n}{\mu_g} \right) K \nabla p^n \right) \cdot \nabla v \, d\Omega = \\
    &\int_{\Omega} \left( \vec{g} K \left( \frac{1 - S_g^n}{\mu_w} \rho_w + \frac{S_g^n}{\mu_g} \rho_g \right) \right) \cdot \nabla v \, d\Omega +
    \int_{\Omega} (q_w^n + q_g^n) v \, d\Omega
\end{aligned}
\end{equation}

\section{Discretization with B-spline basis functions}

Now, we can move on to the numerical implementation of these equations using B-spline basis functions within the IGA-ADS framework.
By using the same logic as for the weak form derivation in the previous section, the equation \eqref{eq:6} can be reformulated as a variational formulation for computing $L^2$ products with test functions $v$:
\begin{eqnarray}{\label{eq:10}}
\begin{aligned}
    &(S_{n+1}\textrm{,}v)_{L2} = ({S_n} + \tau {\cal L} S_{n}+\tau f_n\textrm{,}v)_{L2}; \\
    &{\cal L} S_{n} = \phi^{-1} \nabla \cdot \left( \frac{S_g^n}{\mu_g} K \left( \nabla p^n - \rho_g \vec{g} \right) \right);\\
    &f_n = \phi^{-1} q_g^n
\end{aligned}
\end{eqnarray}
The saturation functions $S_n$ and $S_{n+1}$ and the test function $v$ are then approximated as linear combinations of the B-spline basis functions of order $r$ and continuity $C^{r-1}$:
\begin{eqnarray}{\label{eq:11}}
\begin{aligned}
    &S_n \approx \sum_{i,j} S_n^{i,j} B_{i,r}^x B_{j,r}^y; \\
    &S_{n+1} \approx \sum_{i,j} S_{n+1}^{i,j} B_{i,r}^x B_{j,r}^y; \\
    &v = B_{k,r}^x B_{l,r}^y
\end{aligned}
\end{eqnarray}

Hence, we end up with a system of linear equations (as isogeometric $L2$ projections) for the saturation functions at each time step:
\begin{equation}{\label{eq:12}}
\begin{aligned}
    & \sum_{i\textrm{,}j} S_{n+1}^{i\textrm{,}j} (B^x_{i\textrm{,}r}(x)  B^y_{j\textrm{,}r}(y)\textrm{,}B^x_{k\textrm{,}r}(x)  B^y_{l\textrm{,}r}(y))_{L2} = \\
    & \left({\sum_{i\textrm{,}j} S_n^{i\textrm{,}j} ( B^x_{i\textrm{,}r}(x)  B^y_{j\textrm{,}r}(y)} + \tau \sum_{i\textrm{,}j} S_n^{i\textrm{,}j} {\cal L}(B^x_{i\textrm{,}r}(x)  B^y_{j\textrm{,}r}(y))+f_n , B^x_{k\textrm{,}r}B^y_{l\textrm{,}r}\right)_{L2}, \forall k\textrm{,}l
\end{aligned}
\end{equation} 
We can express it in a matrix form as follows:
\begin{equation}{\label{eq:13}}
\begin{bmatrix}
    \Prod{B^x_{1\textrm{,}r}B^y_{1\textrm{,}r}}{B^x_{1\textrm{,}r}B^y_{1\textrm{,}r}} & \cdots &  \Prod{B^x_{1\textrm{,}r}B^y_{1\textrm{,}r}}{B^x_{N_x\textrm{,}r}B^y_{N_y\textrm{,}r}} \\
    \Prod{B^x_{2\textrm{,}r}B^y_{1\textrm{,}r}}{B^x_{1\textrm{,}r}B^y_{1\textrm{,}r}} & \cdots &  \Prod{B^x_{2\textrm{,}r}B^y_{1\textrm{,}r}}{B^x_{N_x\textrm{,}r}B^y_{N_y\textrm{,}r}} \\   
    \vdots & \vdots & \vdots &  \vdots \\
    \Prod{B^x_{N_x\textrm{,}r}B^y_{N_y\textrm{,}r}}{B^x_{1\textrm{,}r}B^y_{1\textrm{,}r}} &\cdots &  \Prod{B^x_{N_x\textrm{,}r}B^y_{N_y\textrm{,}r}}{B^x_{N_x\textrm{,}r}B^y_{N_y\textrm{,}r}} \\
\end{bmatrix}
\begin{bmatrix}
    S^{1,1}_{n+1} \\ S^{2,1}_{n+1} \\ \vdots \\ S^{N_x,N_y}_{n+1}\\
\end{bmatrix}  \nonumber
\end{equation}
\begin{equation} 
=\begin{bmatrix}
{\left(\sum_{i\textrm{,}j} S_n^{i\textrm{,}j} ( B^x_{i\textrm{,}r}(x)  B^y_{j\textrm{,}r}(y) + \tau \sum_{i\textrm{,}j} S_n^{i\textrm{,}j} L(B^x_{i\textrm{,}r}(x)  B^y_{j\textrm{,}r}(y))+f_n \textrm{,}B^x_{1\textrm{,}r}B^y_{1\textrm{,}r}\right)_{L2}}  \\
{\left(\sum_{i\textrm{,}j} S_n^{i\textrm{,}j} ( B^x_{i\textrm{,}r}(x)  B^y_{j\textrm{,}r}(y) + \tau \sum_{i\textrm{,}j} S_n^{i\textrm{,}j} L(B^x_{i\textrm{,}r}(x)  B^y_{j\textrm{,}r}(y))+f_n \textrm{,}B^x_{2\textrm{,}r}B^y_{1\textrm{,}r}\right)_{L2}}  \\
\vdots \\
{\left(\sum_{i\textrm{,}j} S_n^{i\textrm{,}j} ( B^x_{i\textrm{,}r}(x)  B^y_{j\textrm{,}r}(y) + \tau \sum_{i\textrm{,}j} S_n^{i\textrm{,}j} L(B^x_{i\textrm{,}r}(x)  B^y_{j\textrm{,}r}(y))+f_n \textrm{,}B^x_{N_x\textrm{,}r}B^y_{N_y\textrm{,}r}\right)_{L2}}  \\
\end{bmatrix}
\end{equation}

The Gram matrix build with B-splines over regular  2D domain $\Omega = \Omega_x \times \Omega_y$ has a Kronecker product structure
\begin{equation}{\label{eq:14}}
\begin{aligned}
    \mathcal{M}_{ijkl} &=
    \left( {B_{ij}},{B_{kl}} \right)_{L^2} =
    \int_\Omega B_{ij}B_{kl}\,\mbox{d}\Omega 
    =\int_\Omega B^x_i(x) B^y_j(y) B^x_k(x) B^y_l(y) \,\mbox{d}\Omega \\
    &= \int_\Omega (B_i B_k)(x)\,(B_j B_l)(y)\,\mbox{d}\Omega  
    = \left(\int_{\Omega_x} B_i B_k \,\mbox{d}x\right)
    \left(\int_{\Omega_y} B_j B_l \,\mbox{d}y\right) 
    = \mathcal{M}^x_{ik} \mathcal{M}^y_{jl}
\end{aligned}
\end{equation}

In other words
$\mathcal{M} = \mathcal{M}^x \otimes \mathcal{M}^y $ and the inverse $\mathcal{M}^{-1} = {\mathcal{M}^x}^{-1} \otimes {\mathcal{M}^y}^{-1} $.

Since $\mathcal{M}^x$ and $\mathcal{M}^y$ are multi-diagonal due to local B-splines support, they can be inverted in a linear cost using Thomas algorithm.

The pressure equation \eqref{eq:9} can be implemented in a similar way, however, we can only solve it for the pressure gradient $\nabla p$ instead of the pressure $p$ itself if we want to keep the same Kronecker product structure of the resulting system of equations.
This is clearly visible from the variational formulation presented below.
\begin{equation}{\label{eq:15}}
\begin{aligned}
    \left(\left(\frac{1 - S_g^n}{\mu_w} + \frac{S_g^n}{\mu_g} \right) K \nabla p^n, \nabla \cdot v \right)_{L^2} = &\left(\vec{g} K \left( \frac{1 - S_g^n}{\mu_w} \rho_w + \frac{S_g^n}{\mu_g} \rho_g \right), \nabla \cdot v \right)_{L^2} - (q_w^n + q_g^n, v)_{L^2}
\end{aligned}
\end{equation}

However, we can still express pressure and saturation as linear combinations of B-spline basis functions:
\begin{eqnarray}{\label{eq:16}}
    S^n_g(x,y) = {\sum_{i\textrm{,}j} S_n^{i\textrm{,}j}  B^x_{i\textrm{,}r}(x)  B^y_{j\textrm{,}r}(y)}; \\
    p^n(x,y) = {\sum_{i\textrm{,}j} p_n^{i\textrm{,}j} B^x_{i\textrm{,}r}(x)  B^y_{j\textrm{,}r}(y)}
\end{eqnarray}

This discretization allows us to factor out the pressure values ($p_n^{i\textrm{,}j}$) as control points of the B-spline basis functions, and apply the gradient operator to the basis functions instead of the pressure itself.
Then, we can define a generalized system of equations for pressure $\cal{S}\cal{P}= \cal{F}$. 
\begin{equation}{\label{eq:18}}
\begin{aligned}
{\cal S}_{\left(i,j;k,l\right)} = \\ \left\{ 
\begin{array}{l}
\left( \left(\frac{1 - {\sum_{i\textrm{,}j} S_n^{i\textrm{,}j} ( B^x_{i\textrm{,}r}(x)  B^y_{j\textrm{,}r}(y))}
}{\mu_w} + \frac{{\sum_{i\textrm{,}j} S_n^{i\textrm{,}j} ( B^x_{i\textrm{,}r}(x)  B^y_{j\textrm{,}r}(y))}
}{\mu_g} \right) K(x,y) \right. \\ 
\left. \nabla \cdot \left( B^x_{i\textrm{,}r}(x)  B^y_{j\textrm{,}r}(y)\right)
, \nabla \cdot \left( B^x_{k\textrm{,}r}(x)  B^y_{l\textrm{,}r}(y) \right)
\right)_{L^2}
\end{array}
\right\}_{\left(i,j;k,l\right)} 
\end{aligned}
\end{equation}

\begin{equation}{\label{eq:19}}
{\cal P} = \begin{bmatrix} p^{1,1}_n \\ \cdots \\ p^{N_x,N_y}_n \end{bmatrix}
\end{equation}

\begin{equation}{\label{eq:20}}
\begin{aligned}
{\cal F}_{\left(k,l\right)} =  \\
\left\{ 
\begin{array}{l}
\left(\vec{g} K(x,y) \left( \frac{1 - {\sum_{i\textrm{,}j} S_n^{i\textrm{,}j} ( B^x_{i\textrm{,}r}(x)  B^y_{j\textrm{,}r}(y))}}{\mu_w} \rho_w + \frac{{\sum_{i\textrm{,}j} S_n^{i\textrm{,}j} ( B^x_{i\textrm{,}r}(x)  B^y_{j\textrm{,}r}(y))}}{\mu_g} \rho_g \right), \right. \\
\left. \nabla \cdot B^x_{k\textrm{,}r}(x)  B^y_{l\textrm{,}r}(y)
 \right)_{L^2} - \left(q_w^n + q_g^n, 
B^x_{k\textrm{,}r}(x)  B^y_{l\textrm{,}r}(y)
\right)_{L^2}
\end{array}
\right\}_{\left(k,l\right)} 
\end{aligned}
\end{equation}

This system can be solved using the MUMPS solver \cite{MUMPS:1,MUMPS:2}, which has an interface integrated into the IGA-ADS software. However, the computational cost of the MUMPS solver applied for the pressure updates is ${\cal O}(N^{1.5}p^2)$ \cite{COLLIER2012353} for the two-dimensional problems, and it dominates the linear ${\cal O}(N)$ computational cost of the IGA-ADS solver applied for saturation updates. Thus, we will replace the MUMPS solver with the CRVPINN solver \cite{los_collocation-based_2025}.

\section{Collocation-based Robust Variational Physics Informed Neural Network solver for pressure}

In this paper, we want to replace the calls to the MUMPS direct solver with the Collocation-based Robust Variational PINN (CRVPINN) solver.
We will introduce the main concepts of the CRVPINN formulation used for the pressure computation of the CO2 sequestration problem. For a more detailed description of the method, along with a library of practical examples, we refer to \cite{los_collocation-based_2025}.

We focus on the pressure computation in the~$CO_2$ sequestration problem
\begin{equation}
\nabla\cdot \left(
\underbrace{
\left( \frac{1-S_g}{\mu_w} + \frac{S_g}{\mu_g} \right)
K}_{\alpha}
\nabla p^n
\right)
=
\underbrace{ \nabla\cdot \left(
\vec{g} K\left( 
    \frac{1-S_g}{\mu_w}\rho_w +
    \frac{S_g}{\mu_g}\rho_g \right)
\right) - (q_w + q_g)
}_{f}
\end{equation}

This problem is equivalent to the Poisson problem with a non-constant diffusion coefficient~$\alpha$
\begin{equation}
    \nabla \cdot \left(\alpha \nabla p^n\right) = f
\end{equation}

We define a discrete computational domain: 
\begin{equation}{\label{eq:21}}
    \Omega_h := \{ (ih,jh) \in (0,1)^2: 0 \leq i \leq N, 0 \leq j \leq N\},
\end{equation}
where $h = 1/N$ is the grid spacing, and $N$ is the number of collocation points in each direction.
The discrete gradients are defined as follows:
\begin{equation}{\label{eq:22}}
\begin{aligned}
\nabla_{+} u_{i,j} := \left(\nabla_{x+} u_{i,j},\nabla_{y+} u_{i,j}\right) := 
\left(\frac{u_{i+1,j}-u_{i,j}}{h},\frac{u_{i,j+1}-u_{i,j}}{h}\right),\\
\nabla_{-} u_{i,j} := \left(\nabla_{x-} u_{i,j},\nabla_{y-} u_{i,j}\right) := 
\left(\frac{u_{i,j}-u_{i-1,j}}{h},\frac{u_{i,j}-u_{i,j-1}}{h}\right),
\end{aligned}
\end{equation}

\newcommand{\ProdA}[3]{\left({#2}, {#3}\right)_{#1}}
\newcommand{\Prodh}[2]{\ProdA{h}{#1}{#2}}
\newcommand{\Proddh}[2]{\ProdA{\nabla, h}{#1}{#2}}

\newcommand{\Norm}[2]{\left\|{#2}\right\|_{#1}}
\newcommand{\Normh}[1]{\Norm{h}{#1}}
\newcommand{\Normdh}[1]{\Norm{\nabla, h}{#1}}

\newcommand{\Dp}{\nabla_{+}}
\newcommand{\Dm}{\nabla_{-}}
\newcommand{\Dxp}{\nabla_{x+}}
\newcommand{\Dxm}{\nabla_{x-}}
\newcommand{\Dyp}{\nabla_{y+}}
\newcommand{\Dym}{\nabla_{y-}}

\newcommand{\Dh}{D_h}
\newcommand{\Dzh}{D_{0,h}}
\newcommand{\Dth}{\widetilde{D}_h}

where $u_{i,j} := u(ih,jh)$.
We introduce the equivalence of the continuous  $L^2$-norm and scalar product:
\begin{equation}
    \Prodh{u}{v} = h^2 \sum_{p \in \Omega_h} u(p)v(p),
    \qquad
    \Normh{u}^2 = \Prodh{u}{u}
\end{equation}

We also introduce a discrete equivalence  of~$H^1_0(\Omega)$ space 
\begin{equation}
    \Dzh = \left\{ u \in \Dh : \left.u\right|_{\partial\Omega_h} = 0 \right\},
    \quad
    \partial\Omega_h = \partial\Omega \cap \Omega_h
\end{equation}

as well as the equivalence of the continuous~$H^1_0$-norm and scalar product:
\begin{equation}
\begin{aligned}
    \Proddh{u}{v} &= \Prodh{\Dxp u}{\Dxp v} + \Prodh{\Dyp u}{\Dyp v}, \quad
    \Normdh{u}^2 &= \Proddh{u}{u}
\end{aligned}
\end{equation}

We also introduce discrete Kronecker delta test functions

\begin{equation}{\label{eq:23}}
\begin{aligned}
    v &= \delta_{i,j} := \displaystyle{\left\{
        \begin{aligned}
        &1 \quad \textrm{ if }x=x_{i,j}, \\
        &0 \quad \textrm{ if }x\neq x_{i,j}.
        \end{aligned} \right.}
%    \alpha_{i,j} &= \left(\frac{1 - S_g^n(x_{i,j})}{\mu_w} + \frac{S_g^n(x_{i,j})}{\mu_g} \right) K(x_{i,j}); \\
%    \text{RHS}_{i,j} &= \nabla_{+} \cdot \left( \vec{g} K(x_{i,j}) \left( \frac{1 - S_g^n(x_{i,j})}{\mu_w} \rho_w + \frac{S_g^n(x_{i,j})}{\mu_g} \rho_g \right) \right) - (q_w^n(x_{i,j}) + q_g^n(x_{i,j}))
\end{aligned}
\end{equation}
where $x_{i,j}=(ih,jh)$.

%Equipped with these inner products, and using the Kronecker deltas as test functions, we derive a weak discrete formulation
%\begin{equation}
%    \Prodh{\Dp p^n}{\Dp v} = \Prodh{f}{v} \quad \forall v \in \Dzh
%\end{equation}
%
%Assuming~$0 < \gamma \leq \alpha \leq C $, we have
%
%\begin{equation}
%    \gamma \Normdh{u}
%    \leq
%    \sup_{\substack{v \in \Dzh\\ v \neq 0}}
%    \frac{\Prodh{\Dp p^n}{\Dp v}}{\Normdh{v}}
%    \leq
%    C \Normdh{u}
%\end{equation}
%
%that is, inf-sup stability with constants~$\gamma$, $C$

Then, we can refer to the weak formulation of the problem defined in equation \eqref{eq:15}, adapting it to our CRVPINN formulation.
We will use the Kronecker delta function $\delta_{i,j}$ as the test function. Hence, the building blocks for the residual function are defined as follows:
\begin{equation}{\label{eq:23}}
\begin{aligned}
%    v &= \delta_{i,j} := \displaystyle{\left\{
%        \begin{aligned}
%        &1 \quad \textrm{ if }x=x_{i,j}, \\
%        &0 \quad \textrm{ if }x\neq x_{i,j}.
%        \end{aligned} \right.} \\
    \alpha_{i,j} &= \left(\frac{1 - S_g^n(x_{i,j})}{\mu_w} + \frac{S_g^n(x_{i,j})}{\mu_g} \right) K(x_{i,j}); \\
    \text{RHS}_{i,j} &= \nabla_{+} \cdot \left( \vec{g} K(x_{i,j}) \left( \frac{1 - S_g^n(x_{i,j})}{\mu_w} \rho_w + \frac{S_g^n(x_{i,j})}{\mu_g} \rho_g \right) \right) - (q_w^n(x_{i,j}) + q_g^n(x_{i,j}))
\end{aligned}
\end{equation}
And trivially, we replace the pressure function $p$ with the output of the neural network $\text{PINN}(x_{i,j})$.
The discrete weak formulation for the pressure equation is:
\begin{equation}{\label{eq:24}}
    (\alpha_{i,j} \nabla_{+} \text{PINN}_{i,j}, \nabla_{+} \delta_{k,l}) = (\text{RHS}_{i,j}, \delta_{k,l}), \forall (k,l) \in \{0,1,...,N\}^2
\end{equation}
Now, we can define a discrete weak residual function for the pressure equation:
\begin{equation}{\label{eq:25}}
    \text{RES}_{k,l} = \sum_{i,j} (\alpha_{i,j} \nabla_{+} \text{PINN}_{i,j}, \nabla_{+} \delta_{k,l}) - (\text{RHS}_{i,j}, \delta_{k,l})
\end{equation}

We notice that assuming~$0 < \gamma \leq \alpha \leq C $, we have

\begin{equation}
    \gamma \Normdh{u}
    \leq
    \sup_{\substack{v \in \Dzh\\ v \neq 0}}
    \frac{\Prodh{\Dp p^n}{\Dp v}}{\Normdh{v}}
    \leq
    C \Normdh{u}
\end{equation}

that is, we can show inf-sup stability with constants~$\gamma$, $C$, using $\Prodh{\Dp p^n}{\Dp v}=\Proddh{p^n}{v}$. 
To move forward from a discrete weak residual to a discrete robust loss function, we need to define the Gram matrix
constructed with the inner product of Kronecker deltas $ \left(\delta_{ij}, \delta_{kl} \right)_{\nabla,h}$.
\begin{eqnarray}{\label{eq:26}}
    \mathbf{G}_{i,j;k,l}=
    &h^{-2}\displaystyle{\left\{\begin{aligned}
    &\phantom{-}\,\,4\ \quad \textrm{ for }(i,j)=(k,l), \\
    &-1 \quad \textrm{ for }(k,l)\in \{(i+1,j),(i-1,j)\},\\
    &-1 \quad \textrm{ for }(k,l)\in\{(i,j+1),(i,j-1)\}.
    \end{aligned}
    \right.} 
\end{eqnarray}
Finally, our loss function for the pressure equation is:
\begin{equation}{\label{eq:27}}
    \text{LOSS}(\theta) = \text{RES}^T(\theta) \mathbf{G}^{-1} \text{RES}(\theta)
\end{equation}
Instead of direct inversion of the Gram matrix $\mathbf{G}$, we perform LU factorization, as shown below:
\begin{equation}{\label{eq:28}}
\begin{aligned}
    \mathbf{G} &= \mathbf{L} \mathbf{U} \\
    \mathbf{U} z &= \text{RES}(\theta) \\
    \mathbf{L} q &= z \\
    q &= \mathbf{L}^{-1} z = \mathbf{L}^{-1} \mathbf{U}^{-1} \text{RES}(\theta) = \mathbf{G}^{-1} \text{RES}(\theta) \\
    \text{LOSS}(\theta) & = \text{RES}^T(\theta) q
\end{aligned}
\end{equation}
Here, in order to find $z$ and $q$, we only need to perform forward and backward substitution, which have linear computational cost. The LU factorization of the Gram matrix is performed only once at the beginning of the simulation.

\section{IGA-ADS-CRVPINN Algorithm}

\begin{enumerate}
\item Given the permeability map $K$, the porosity map $\phi$, the source terms $q_w$, $q_g$, the gravitational acceleration $\vec{g}$, and the viscosities $\mu_w$, $\mu_g$ of the brine solution and the gas, along with the initial saturation $S^0_g$ expressed as a linear combination of B-splines.
\item Run the pretraining of the CRVPINN for the initial pressure $p^0$. Minimize the loss 
\eqref{eq:27} with the Gram matrix ${\bf G}$ given by \eqref{eq:26}, using the vectorized residual, namely
\begin{equation}
\text{LOSS}(\theta)=\sum_{i=1,...,N}\sum_{j=1,...,N}\sum_{k=1,...,N}\sum_{l=1,...,N}\mathbf{G}_{i,j;k,l}\text{RES}_{i,j}(\theta)\text{RES}_{k,l}(\theta)
\end{equation}
where 
\begin{equation}
    \text{RES}_{k,l} = \sum_{i,j} (\alpha_{i,j} \nabla_{+} \text{PINN}_{i,j}(\theta), \nabla_{+} \delta_{k,l}) - (\text{RHS}_{i,j}, \delta_{k,l})
\end{equation}
with the right-hand side \eqref{eq:25} including the initial saturation $S^0_g$
\begin{equation}
    \text{RHS}_{i,j} = \nabla_{+} \cdot \left( \vec{g} K(x_{i,j}) \left( \frac{1 - S_g^0(x_{i,j})}{\mu_w} \rho_w + \frac{S_g^0(x_{i,j})}{\mu_g} \rho_g \right) \right) - (q_w^0(x_{i,j}) + q_g^0(x_{i,j}))
\end{equation}
\item Use $n=1$ to obtain $S^n_g=S^0_g$.
\item Compute the $L^2$ projection of the resulting discrete pressure approximated by $\text{PINN}_{i,j}$ into the B-splines representation $p^n$. 
\item Compute the saturation update from \eqref{eq:6}
\begin{equation}
S_g^{n+1} = S_g^n + \phi^{-1} \tau \nabla \cdot \left( \frac{S_g^n}{\mu_g} K \left( \nabla p^n - \rho_g \vec{g} \right) \right) + \phi^{-1} \tau q_g^n
\end{equation}
\item Repeat point 5 ten times using the same pressure approximation.%\footnote{Note. This approximation is considered adequate given that the finite element method solver requires smaller time steps for accurate saturation solutions, while the CRVPINN Pressure solver can approximate the solution without a noticeable drop in accuracy and with a larger time step. This allows us to reduce the computational cost while keeping the same level of quality of the results. }
%\item Train the update of the pressure approximation $p^{n+1}$ using the CRVPINN. 
\item Train the update of the pressure approximation $p^{n+1}$ using the CRVPINN. Minimize the loss 
\eqref{eq:27} with the Gram matrix ${\bf G}$ given by \eqref{eq:26}, using the vectorized residual, 
with the right-hand side \eqref{eq:25} including the actual saturation $S_g^n$ approximation
\begin{equation}
    \text{RHS}_{i,j} = \nabla_{+} \cdot \left( \vec{g} K(x_{i,j}) \left( \frac{1 - S_g^n(x_{i,j})}{\mu_w} \rho_w + \frac{S_g^n(x_{i,j})}{\mu_g} \rho_g \right) \right) - (q_w^n(x_{i,j}) + q_g^n(x_{i,j}))
\end{equation}
\item If the maximum number of iterations has been reached, then STOP. Otherwise, go to 5.
\end{enumerate}

\section{IGA-ADS-MUMPS Algorithm}

\begin{enumerate}
\item Given the permeability map $K$, the porosity map $\phi$, the source terms $q_w$, $q_g$, the gravitational acceleration $\vec{g}$, and the viscosities $\mu_w$, $\mu_g$ of the brine solution and the gas, along with the initial saturation $S^0_g$ expressed as the linear combination of B-splines.
\item Solve the initial pressure configuration 
from \eqref{eq:18}-\eqref{eq:20} using the MUMPS solver
\begin{equation}
\left\{ 
\begin{array}{l}
\left( \left(\frac{1 - {\sum_{i\textrm{,}j} S_0^{i\textrm{,}j} ( B^x_{i\textrm{,}r}(x)  B^y_{j\textrm{,}r}(y))}
}{\mu_w} + \frac{{\sum_{i\textrm{,}j} S_0^{i\textrm{,}j} ( B^x_{i\textrm{,}r}(x)  B^y_{j\textrm{,}r}(y))}
}{\mu_g} \right) K(x,y) \right. \\ 
\left. \nabla \cdot \left( B^x_{i\textrm{,}r}(x)  B^y_{j\textrm{,}r}(y)\right)
, \nabla \cdot \left( B^x_{k\textrm{,}r}(x)  B^y_{l\textrm{,}r}(y) \right)
\right)_{L^2}
\end{array}
\right\}_{\left(i,j;k,l\right)} 
\begin{bmatrix} p^{1,1}_0 \\ \cdots \\ p^{N_x,N_y}_0 \end{bmatrix}= \notag
\end{equation}

\begin{equation}
\begin{aligned}
\left\{ 
\begin{array}{l}
\left(\vec{g} K(x,y) \left( \frac{1 - {\sum_{i\textrm{,}j} S_0^{i\textrm{,}j} ( B^x_{i\textrm{,}r}(x)  B^y_{j\textrm{,}r}(y))}
}{\mu_w} \rho_w + \frac{{\sum_{i\textrm{,}j} S_0^{i\textrm{,}j} ( B^x_{i\textrm{,}r}(x)  B^y_{j\textrm{,}r}(y))}
}{\mu_g} \rho_g \right), \right. \\
\left. \nabla \cdot B^x_{k\textrm{,}r}(x)  B^y_{l\textrm{,}r}(y)
 \right)_{L^2} - \left(q_w^0 + q_g^0, 
B^x_{k\textrm{,}r}(x)  B^y_{l\textrm{,}r}(y)
\right)_{L^2}
\end{array}
\right\}_{\left(k,l\right)} 
\end{aligned}
\end{equation}
\item Define $n=0$ so $S^n_g=S^0_g$ and $p^n=p^0$.
\item Compute the saturation update from \eqref{eq:6}
\begin{equation}
S_g^{n+1} = S_g^n + \phi^{-1} \tau \nabla \cdot \left( \frac{S_g^n}{\mu_g} K \left( \nabla p^n - \rho_g \vec{g} \right) \right) + \phi^{-1} \tau q_g^n
\end{equation}
%\item Repeat point 4 ten times using the same pressure approximation.
\item Solve the pressure update 
from \eqref{eq:18}-\eqref{eq:20} using the MUMPS solver (similar to point 2, but with the actual saturation approximation $S^n_g$). 
\begin{equation}
 \left\{ 
\begin{array}{l}
\left( \left(\frac{1 - {\sum_{i\textrm{,}j} S_n^{i\textrm{,}j} ( B^x_{i\textrm{,}r}(x)  B^y_{j\textrm{,}r}(y))}
}{\mu_w} + \frac{{\sum_{i\textrm{,}j} S_n^{i\textrm{,}j} ( B^x_{i\textrm{,}r}(x)  B^y_{j\textrm{,}r}(y))}
}{\mu_g} \right) K(x,y) \right. \\ 
\left. \nabla \cdot \left( B^x_{i\textrm{,}r}(x)  B^y_{j\textrm{,}r}(y)\right)
, \nabla \cdot \left( B^x_{k\textrm{,}r}(x)  B^y_{l\textrm{,}r}(y) \right)
\right)_{L^2}
\end{array}
\right\}_{\left(i,j;k,l\right)} 
\begin{bmatrix} p^{1,1}_n \\ \cdots \\ p^{N_x,N_y}_n \end{bmatrix} = \notag
\end{equation}

\begin{equation}
\begin{aligned}
\left\{ 
\begin{array}{l}
\left(\vec{g} K(x,y) \left( \frac{1 - {\sum_{i\textrm{,}j} S_n^{i\textrm{,}j} ( B^x_{i\textrm{,}r}(x)  B^y_{j\textrm{,}r}(y))}
}{\mu_w} \rho_w + \frac{{\sum_{i\textrm{,}j} S_n^{i\textrm{,}j} ( B^x_{i\textrm{,}r}(x)  B^y_{j\textrm{,}r}(y))}
}{\mu_g} \rho_g \right), \right. \\
\left. \nabla \cdot B^x_{k\textrm{,}r}(x)  B^y_{l\textrm{,}r}(y)
 \right)_{L^2} - \left(q_w^n + q_g^n, 
B^x_{k\textrm{,}r}(x)  B^y_{l\textrm{,}r}(y)
\right)_{L^2}
\end{array}
\right\}_{\left(k,l\right)} 
\end{aligned}
\end{equation}
\item If the maximum number of iterations has been reached, then STOP. Otherwise, go to 4.
\end{enumerate}

\section{Results}

\subsection{Simulation Setup}

\begin{figure}[htbp]
    \centering
    \begin{subfigure}[b]{0.32\textwidth}
        \centering
        \includegraphics[width=\textwidth]{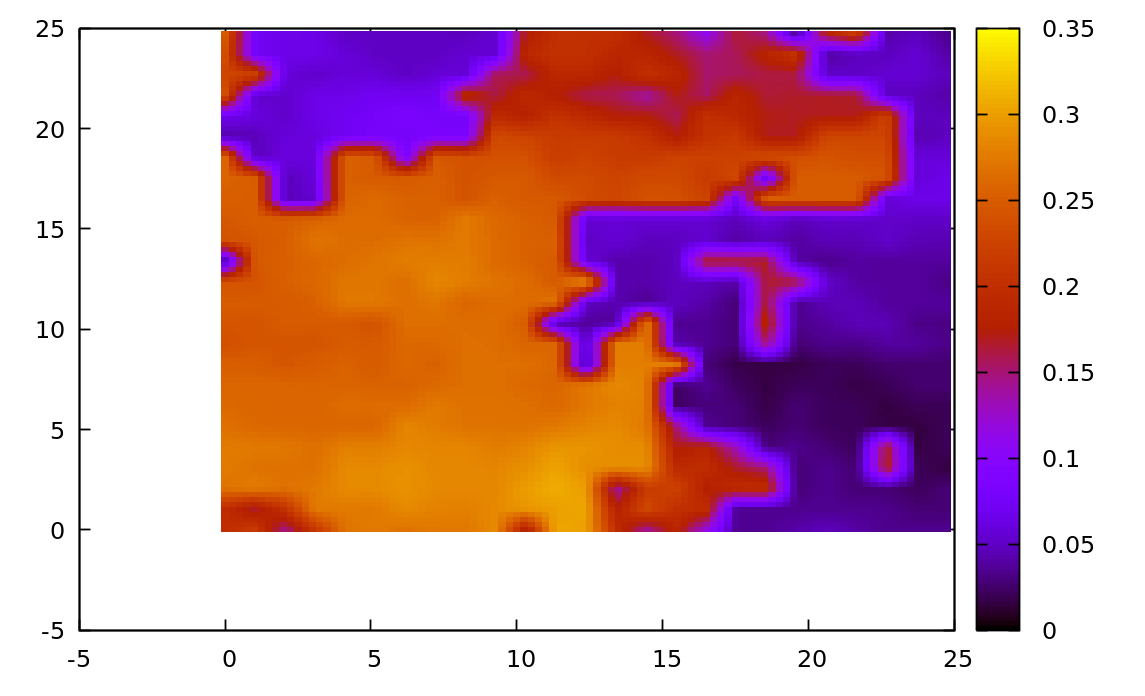}
        \caption{Configuration K1}
        \label{fig:porosity_K1}
    \end{subfigure}
    \hfill
    \begin{subfigure}[b]{0.32\textwidth}
        \centering
        \includegraphics[width=\textwidth]{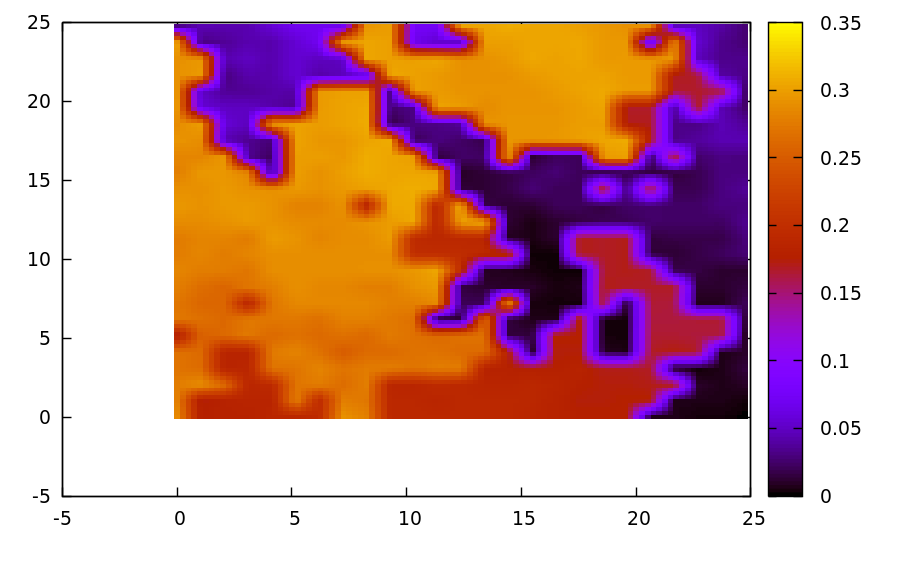}
        \caption{Configuration K2}
        \label{fig:porosity_K2}
    \end{subfigure}
    \hfill
    \begin{subfigure}[b]{0.32\textwidth}
        \centering
        \includegraphics[width=\textwidth]{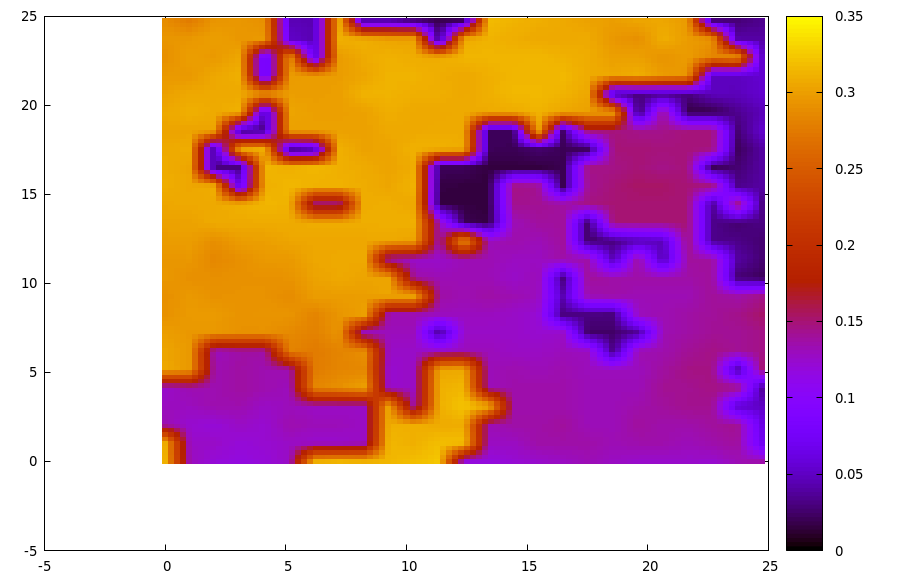}
        \caption{Configuration K3}
        \label{fig:porosity_K3}
    \end{subfigure}
    \caption{Porosity distributions for three different reservoir configurations. The values range from 0 to 1.}
    \label{fig:porosity_maps}
\end{figure}

\begin{figure}[htbp]
    \centering
    \begin{subfigure}[b]{0.35\textwidth}
        \centering
        \includegraphics[width=\textwidth]{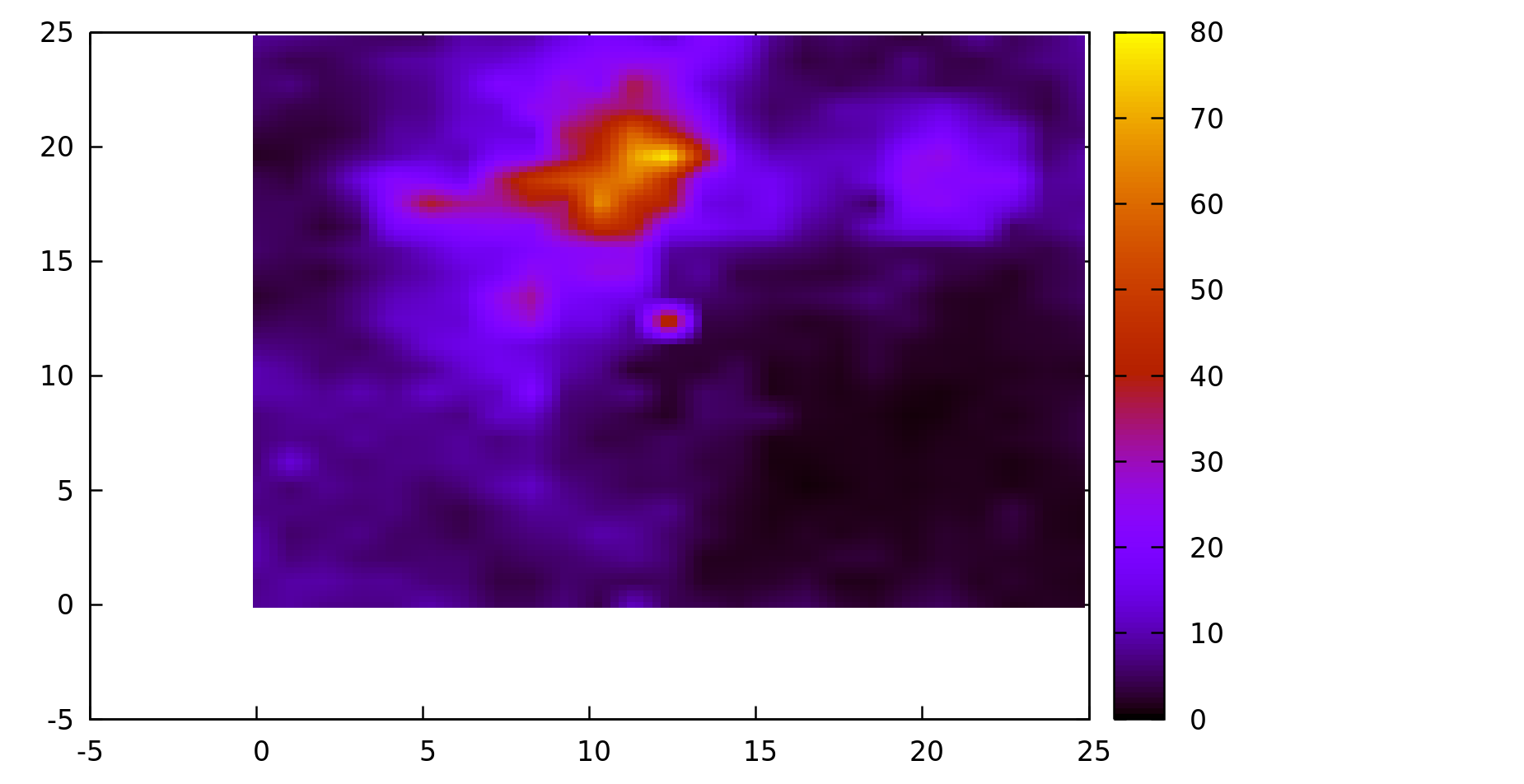}
        \caption{Configuration K1}
        \label{fig:permeability_K1}
    \end{subfigure}
    \hfill
    \begin{subfigure}[b]{0.31\textwidth}
        \centering
        \includegraphics[width=\textwidth]{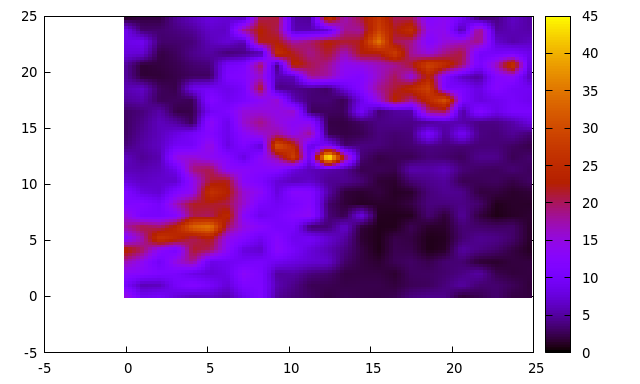}
        \caption{Configuration K2}
        \label{fig:permeability_K2}
    \end{subfigure}
    \hfill
    \begin{subfigure}[b]{0.31\textwidth}
        \centering
        \includegraphics[width=\textwidth]{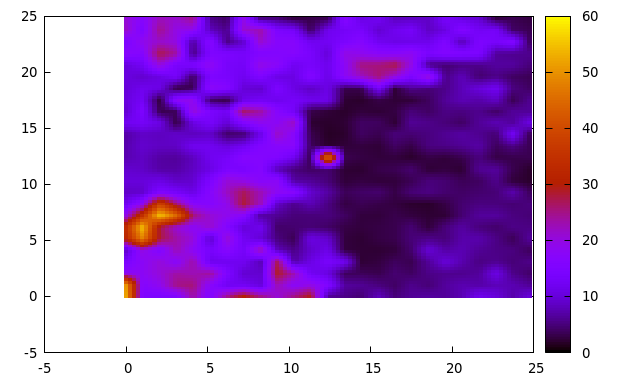}
        \caption{Configuration K3}
        \label{fig:permeability_K3}
    \end{subfigure}
    \caption{Permeability distributions for three different reservoir configurations. The units are in mDarcy, in simulation we use SI units (1 mDarcy $\simeq 10^{-15}$ m$^2$).}
    \label{fig:permeability_maps}
\end{figure}

We present results for three different test configurations (K1, K2, K3) with varying porosity and permeability maps.
The porosity and permeability distributions for each configuration are illustrated in Figures \ref{fig:porosity_maps} and \ref{fig:permeability_maps}, respectively. 
For all our simulations, if we use a permeability map, we use the corresponding porosity map (i.e , K1 uses porosity map 1, K2 uses porosity map 2, etc.).
The computational domain is a 2D rectangular region of size $50 m \times 50 m$.
The physical and numerical parameters used in the simulations are summarized in Table \eqref{tab:parameters}.

The simulations are performed using a uniform mesh with appropriate boundary conditions: no-flux conditions, and pressure boundary conditions on the left, right, top, and bottom boundaries. 
The $CO_2$ injection is modeled through a source term located at the center of the domain, representing an injection well.

\begin{table}[h!]
\centering
    \begin{tabular}{|c|c|}
    \hline
     & SI \\ \hline
    $\rho_g$ & 479 (kg/m$^3$) \\ \hline
    $\rho_w$ & 1045 (kg/m$^3$) \\ \hline
    $g$ & 9.81 (m/s$^2$) \\ \hline
    $\mu_g$ & $3.95 \times 10^{-5}$ (Pa $\cdot$ s) \\ \hline
    $\mu_w$ & $25.35 \times 10^{-5}$ (Pa $\cdot$ s) \\ \hline
    $K_A$ & $1 \times 10^{-15}$ (m$^2$) \\ \hline
%    $K_L$ & $1 \times 10^{-12}$ (m$^2$) & $1 \times 10^{-8}$ (m$^2$) \\ \hline
    \end{tabular}
    \caption{Simulation parameters in SI  units}
\label{tab:parameters}
\end{table}

\subsection{Numerical experiments with IGA-ADS-MUMPS solver}

\begin{figure}[htbp]
    \centering
    \begin{subfigure}[b]{0.9\textwidth}
        \centering
        \includegraphics[height=0.2\textheight]{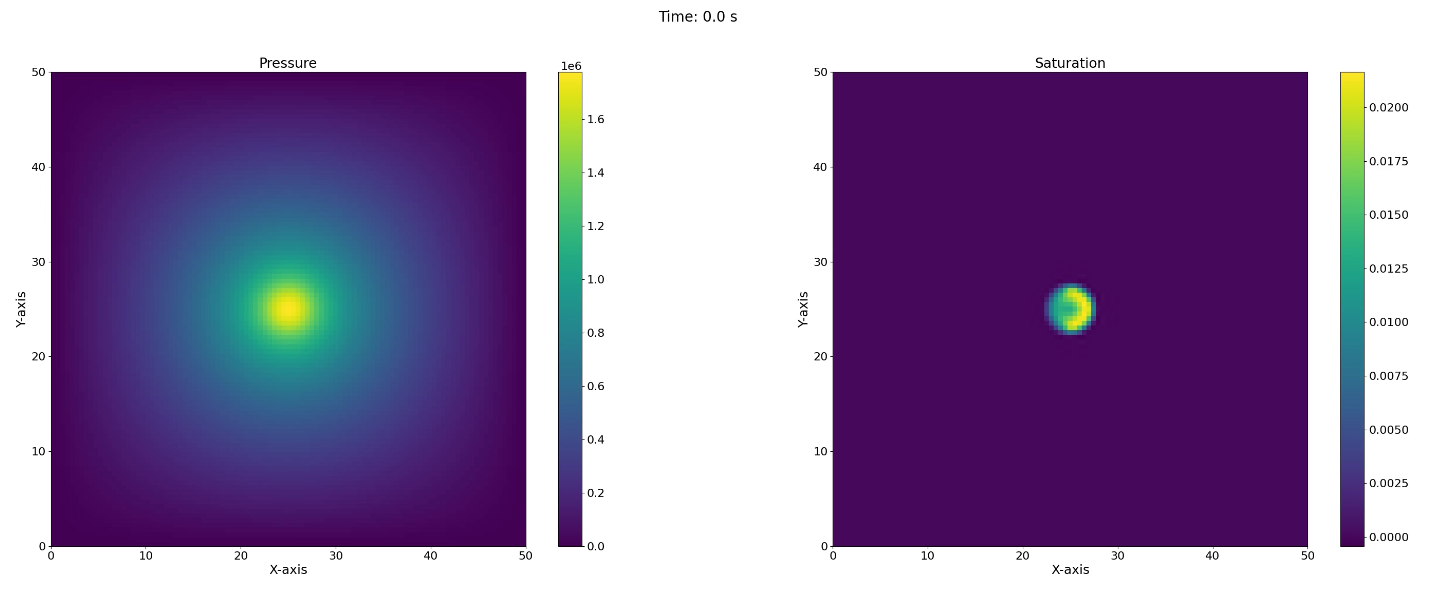}
        \caption{t = 0 s}
        \label{fig:K1_uniform_t0}
    \end{subfigure}
    \hfill
    \begin{subfigure}[b]{0.9\textwidth}
        \centering
        \includegraphics[height=0.2\textheight]{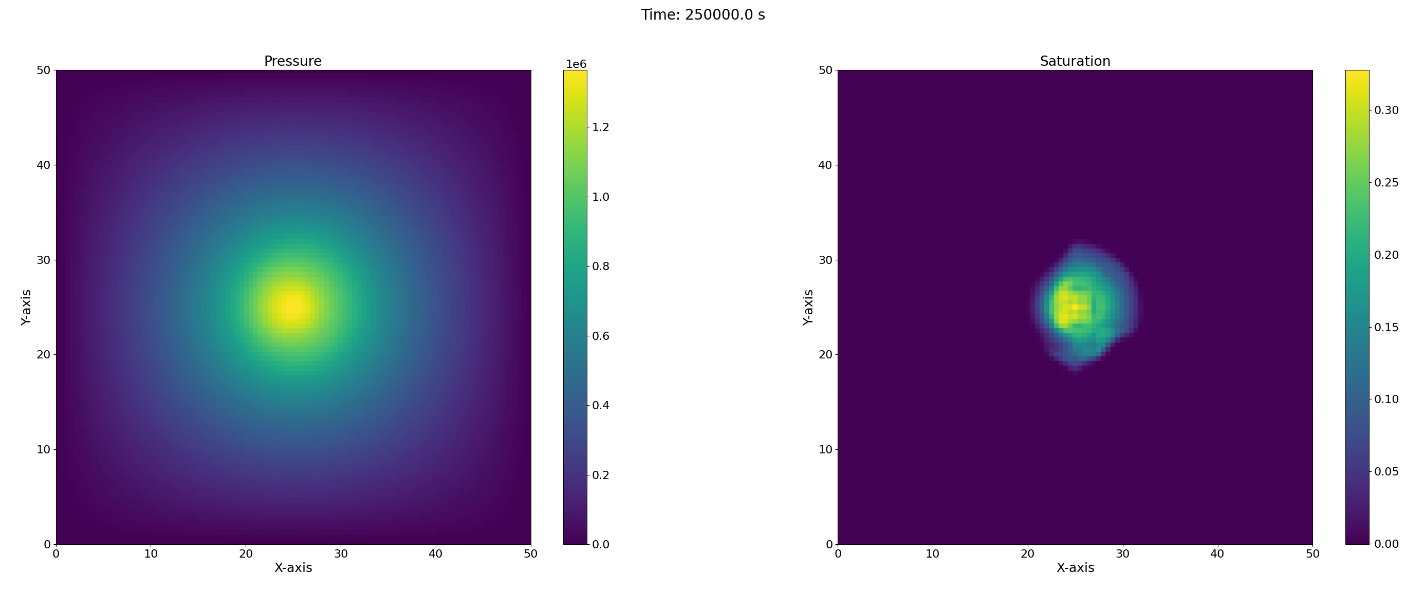}
        \caption{t = 250,000 s (= 50 time steps)}
        \label{fig:K1_uniform_t100}
    \end{subfigure}
    \\
    \begin{subfigure}[b]{0.9\textwidth}
        \centering
        \includegraphics[height=0.2\textheight]{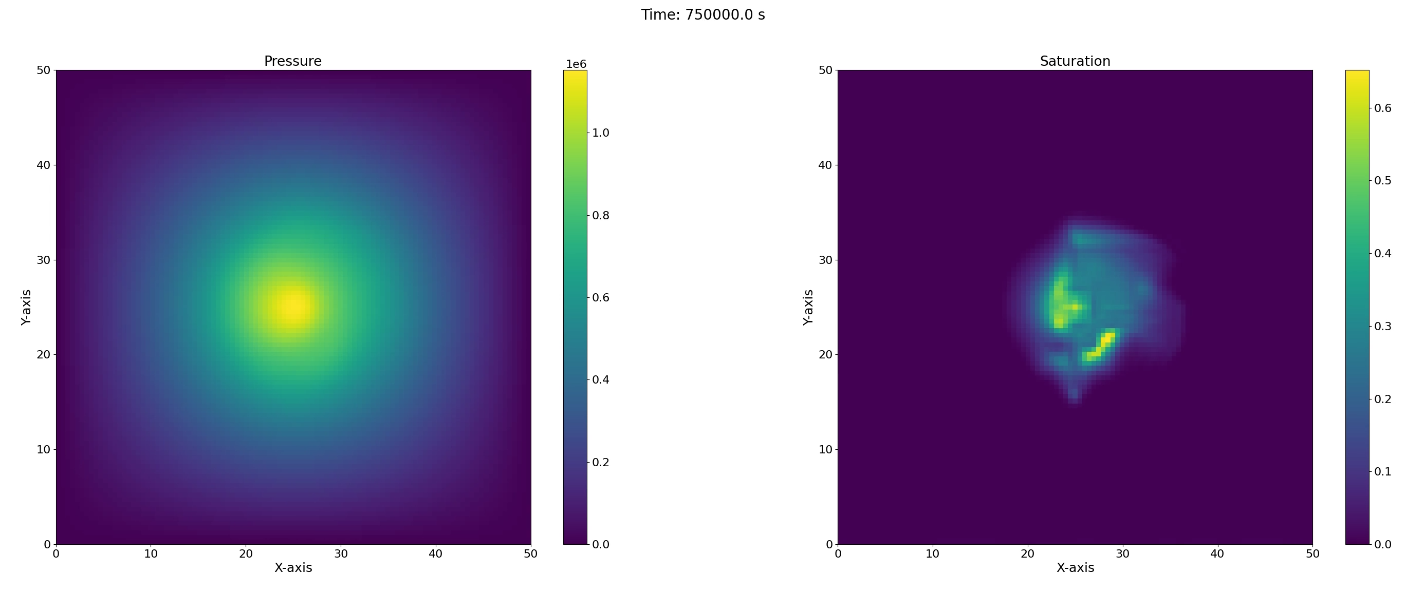}
        \caption{t = 750,000 s (= 150 time steps)}
        \label{fig:K1_uniform_t200}
    \end{subfigure}
    \hfill
    \begin{subfigure}[b]{0.9\textwidth}
        \centering
        \includegraphics[height=0.2\textheight]{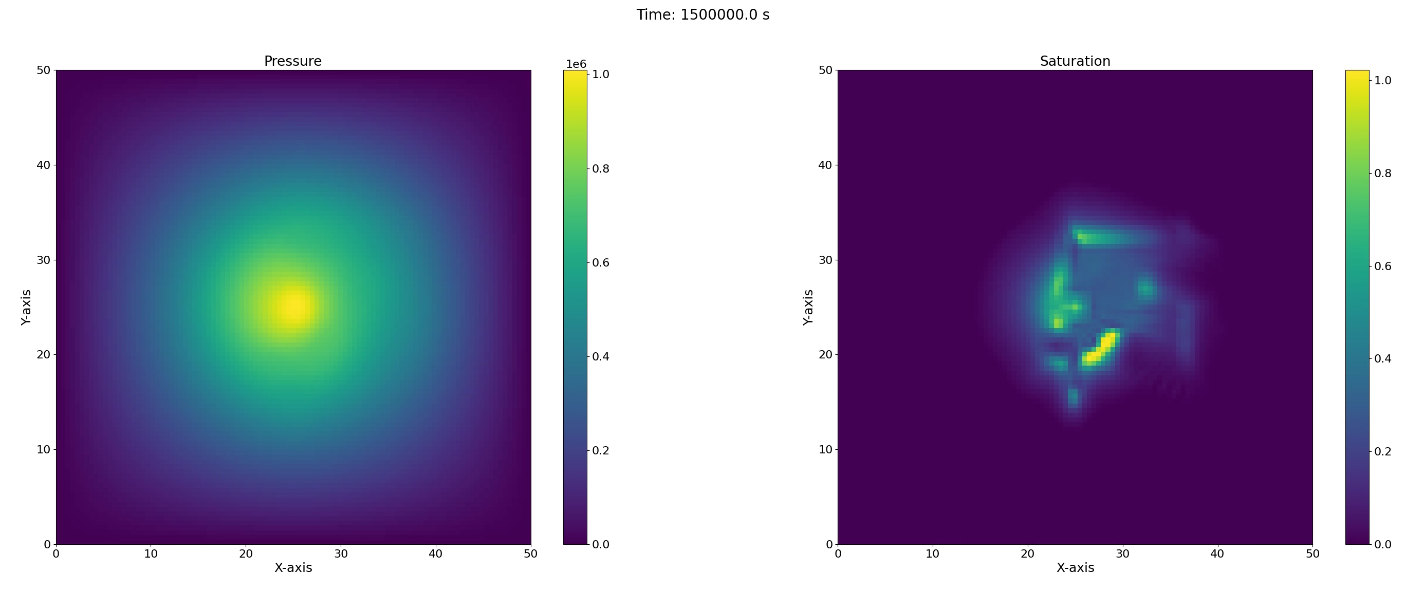}
        \caption{t = 1,250,000 s (= 300 time steps)} 
        \label{fig:K1_uniform_t400}
    \end{subfigure}
    \caption{CO2 saturation evolution for configuration K1 with uniform permeability at different time steps.}
    \label{fig:results_uniformK1}
\end{figure}

\begin{figure}[htbp]
    \centering
    \begin{subfigure}[b]{0.9\textwidth}
        \centering
        \includegraphics[height=0.2\textheight]{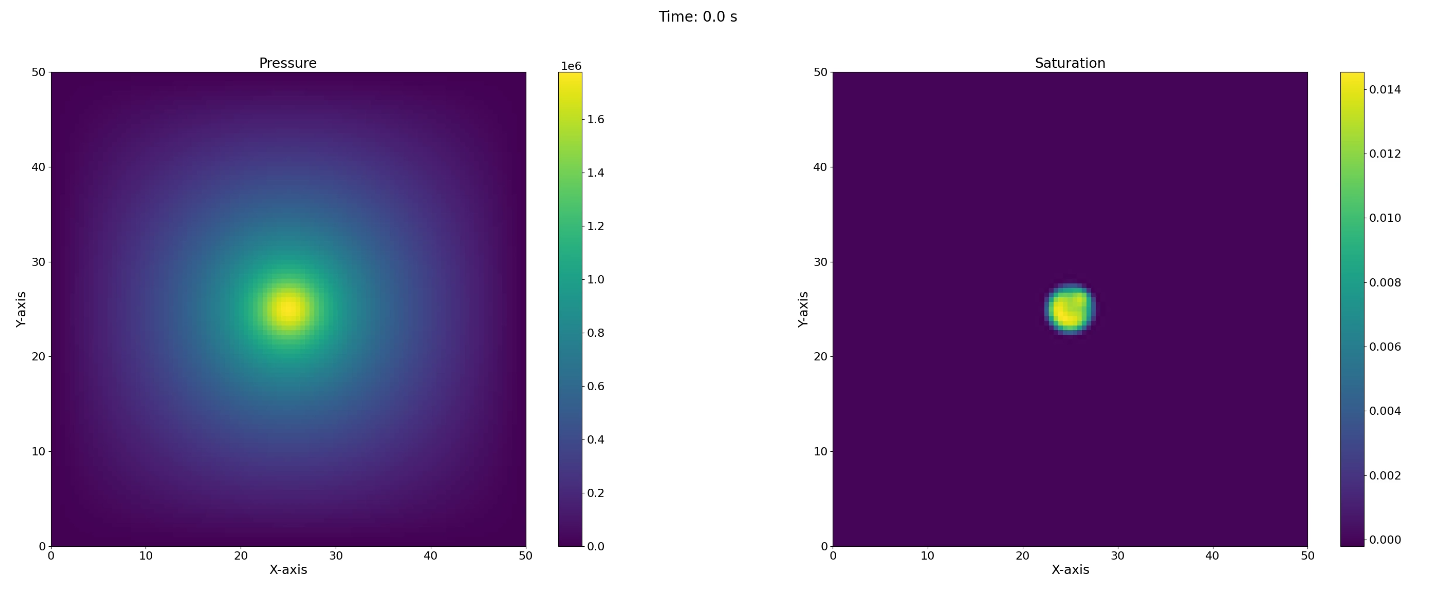}
        \caption{t = 0 s}
        \label{fig:K2_uniform_t0}
    \end{subfigure}
    \hfill
    \begin{subfigure}[b]{0.9\textwidth}
        \centering
        \includegraphics[height=0.2\textheight]{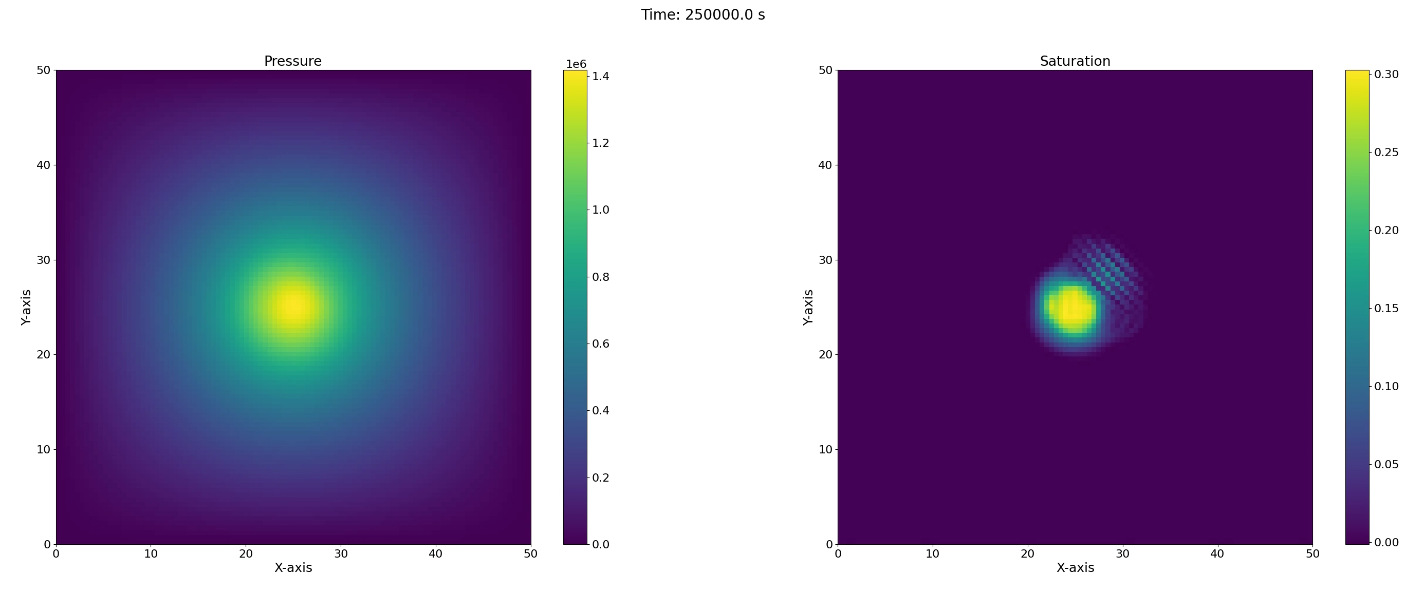}
        \caption{t = 250,000 s (= 50 time steps)}
        \label{fig:K2_uniform_t100}
    \end{subfigure}
    \\
    \begin{subfigure}[b]{0.9\textwidth}
        \centering
        \includegraphics[height=0.2\textheight]{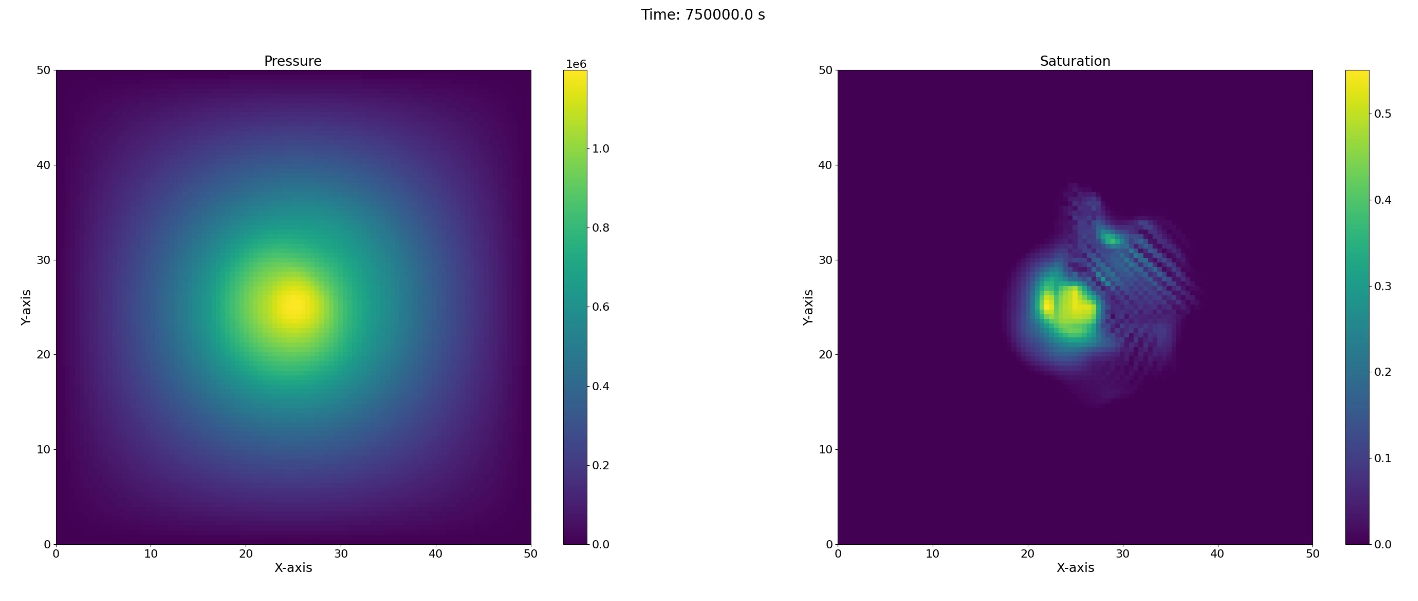}
        \caption{t = 750,000 s (= 150 time steps)}
        \label{fig:K2_uniform_t200}
    \end{subfigure}
    \hfill
    \begin{subfigure}[b]{0.9\textwidth}
        \centering
        \includegraphics[height=0.2\textheight]{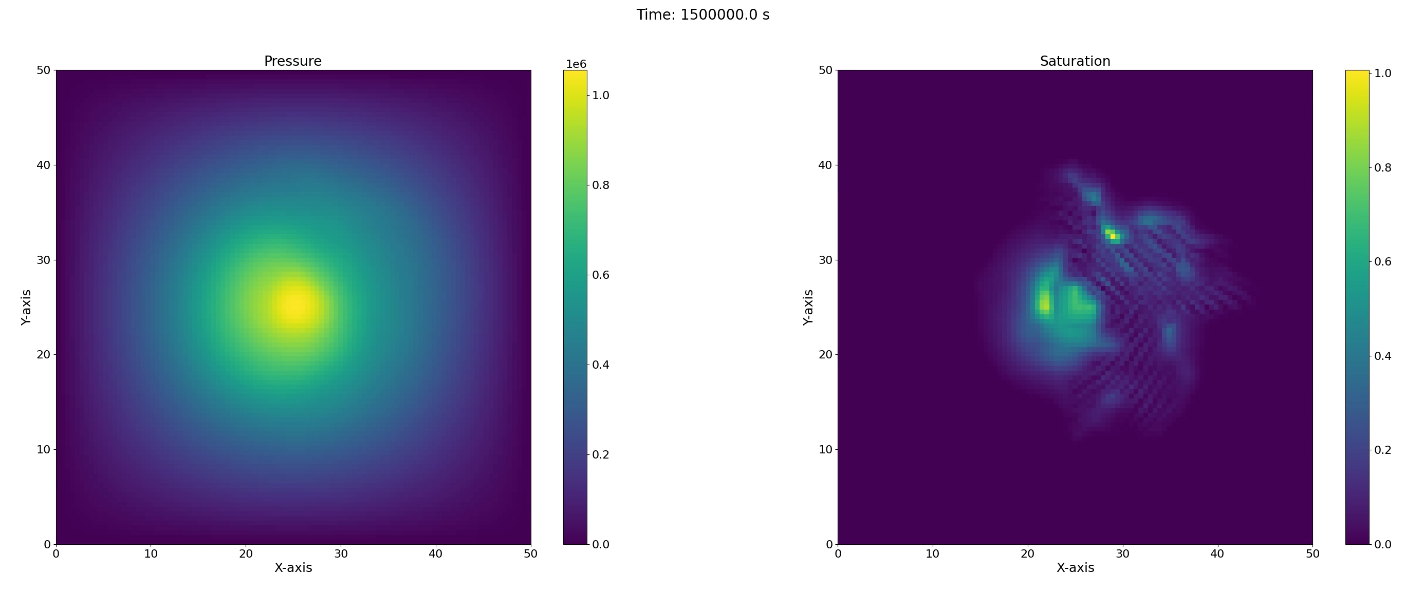}
        \caption{t = 1,250,000 s (= 300 time steps)}
        \label{fig:K2_uniform_t400}
    \end{subfigure}
    \caption{CO2 saturation evolution for configuration K2 with uniform permeability at different time steps.}
    \label{fig:results_uniformK2}
\end{figure}

\begin{figure}[htbp]
    \centering
    \begin{subfigure}[b]{0.9\textwidth}
        \centering
        \includegraphics[height=0.2\textheight]{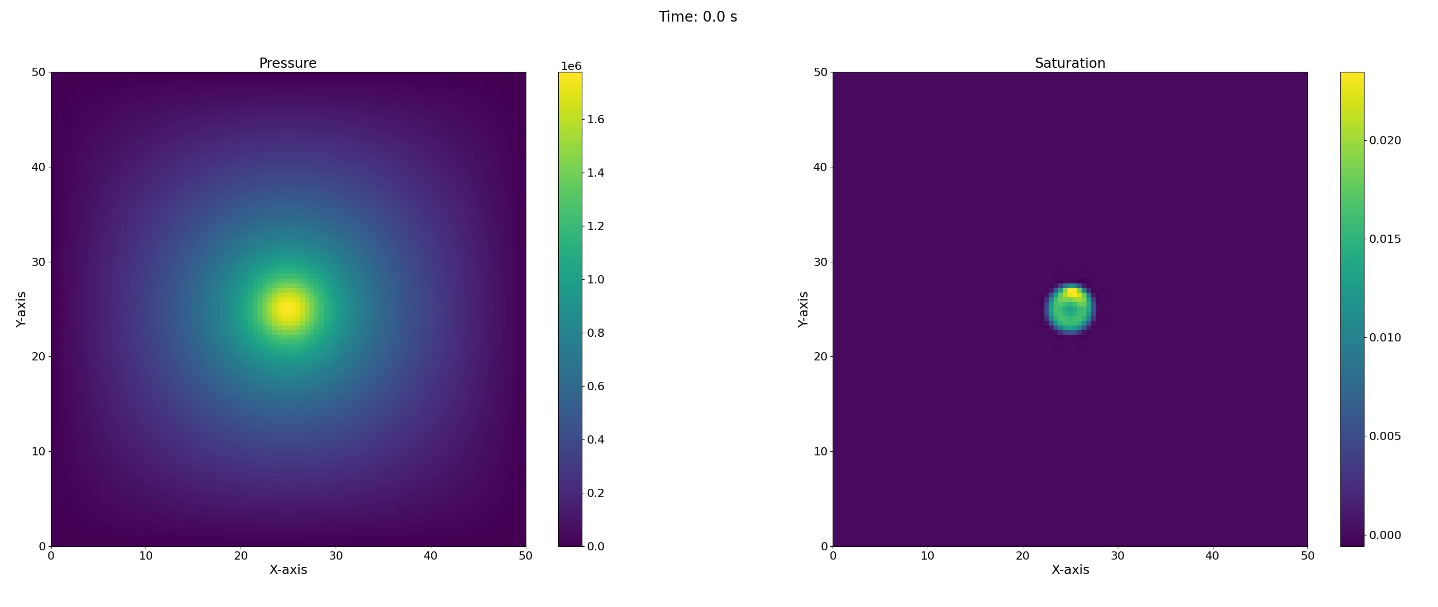}
        \caption{t = 0 s}
        \label{fig:K3_uniform_t0}
    \end{subfigure}
    \hfill
    \begin{subfigure}[b]{0.9\textwidth}
        \centering
        \includegraphics[height=0.2\textheight]{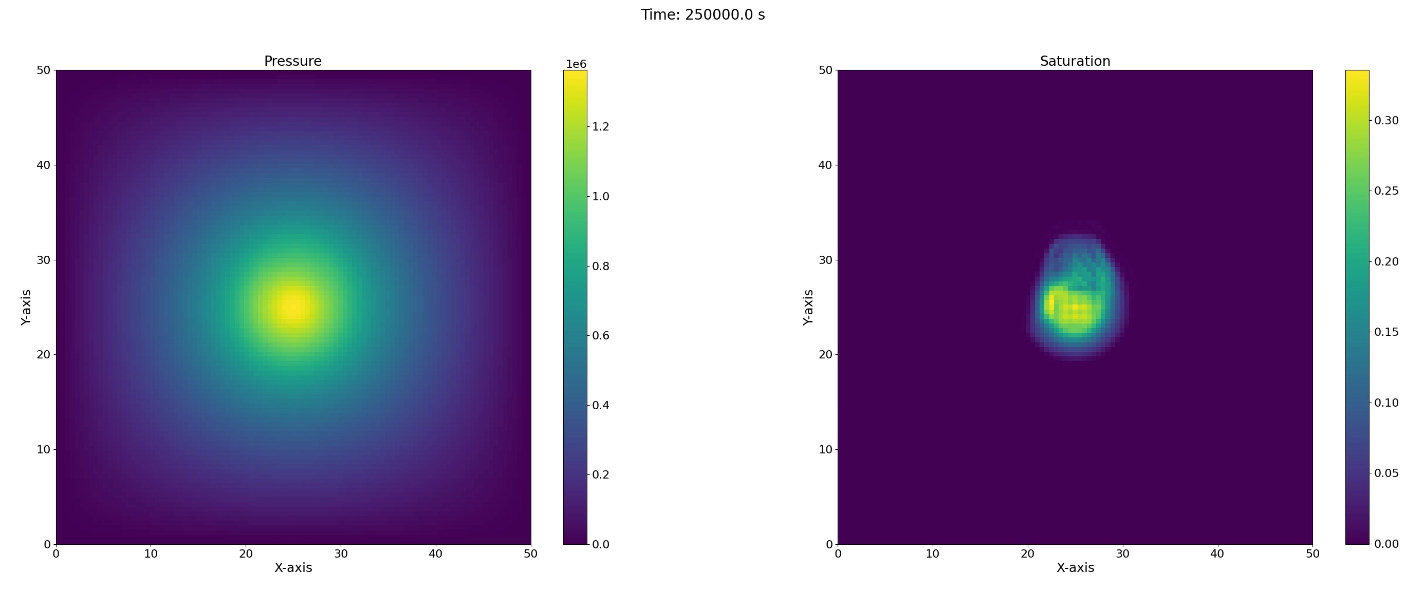}
        \caption{t = 250,000 s (= 50 time steps)}
        \label{fig:K3_uniform_t100}
    \end{subfigure}
    \\
    \begin{subfigure}[b]{0.9\textwidth}
        \centering
        \includegraphics[height=0.2\textheight]{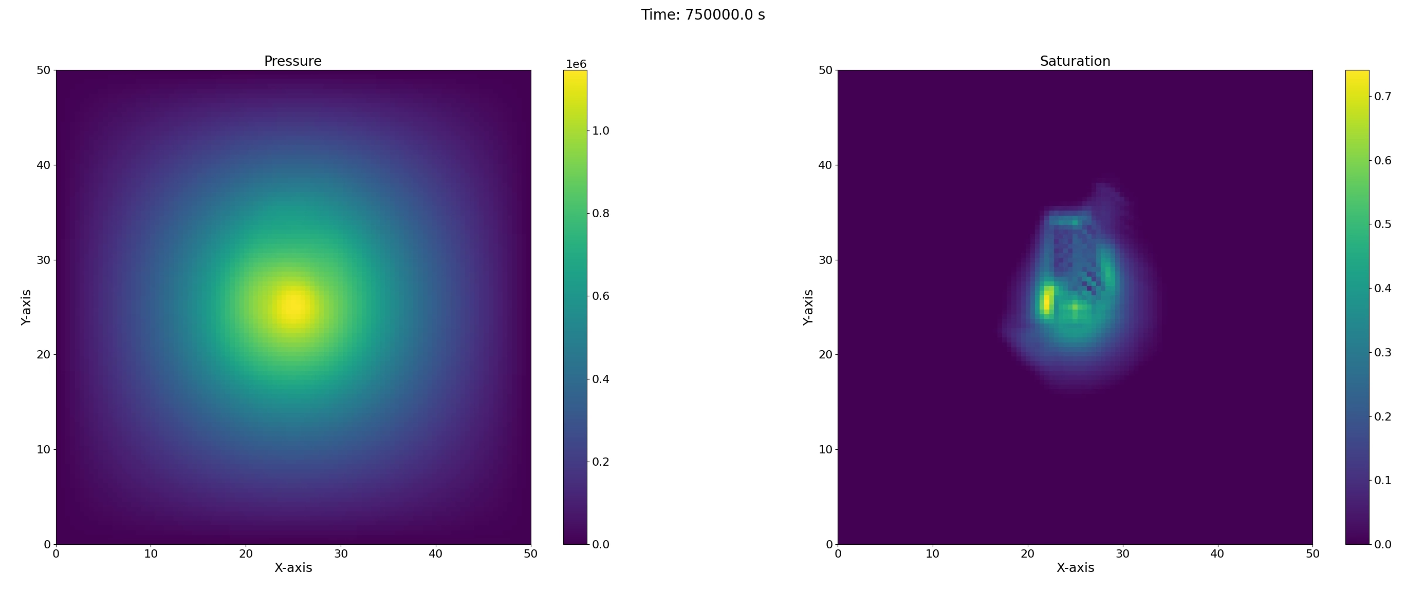}
        \caption{t = 750,000 s (= 150 time steps)}
        \label{fig:K3_uniform_t200}
    \end{subfigure}
    \hfill
    \begin{subfigure}[b]{0.9\textwidth}
        \centering
        \includegraphics[height=0.2\textheight]{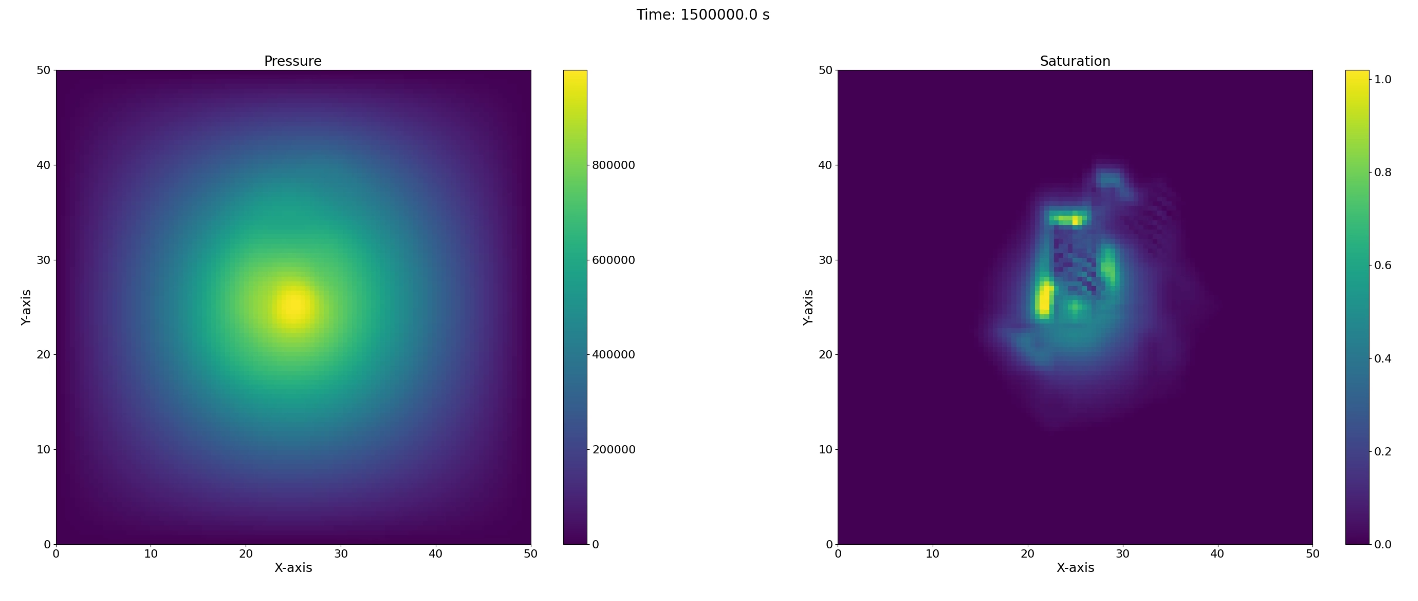}
        \caption{t = 1,250,000 s (= 300 time steps)}
        \label{fig:K3_uniform_t400}
    \end{subfigure}
    \caption{CO2 saturation evolution for configuration K3 with uniform permeability at different time steps.}
    \label{fig:results_uniformK3}
\end{figure}

First, we use the IGA-ADS solver to solve the equation \eqref{eq:6} iteratively for each time step, and the MUMPS solver to solve the equation \eqref{eq:5}.
We completed simulations for the configurations K1, K2, and K3 using only a porosity map and both porosity and permeability maps.
For the porosity-only case, we set a constant permeability value of $K = 1 \times 10^{-15}$ m$^2$ throughout the domain. 
The simulation parameters are set as follows: time step size $\tau = 5000$ seconds, number of time steps = 500, and the source term is defined as a circle of radius 3 m at the center of the domain with a strength of $q_g = 1 \times 10^{-6}$ (dimensionless units of $S_g$) per iteration.
The results of the simulations are presented in Figures \ref{fig:results_uniformK1}, \ref{fig:results_uniformK2}, \ref{fig:results_uniformK3}.
We employed 100 times 100 intervals for the IGA-ADS solver, using quadratic B-splines.

\begin{figure}[htbp]
    \centering
    \begin{subfigure}[b]{0.9\textwidth}
        \centering
        \includegraphics[height=0.2\textheight]{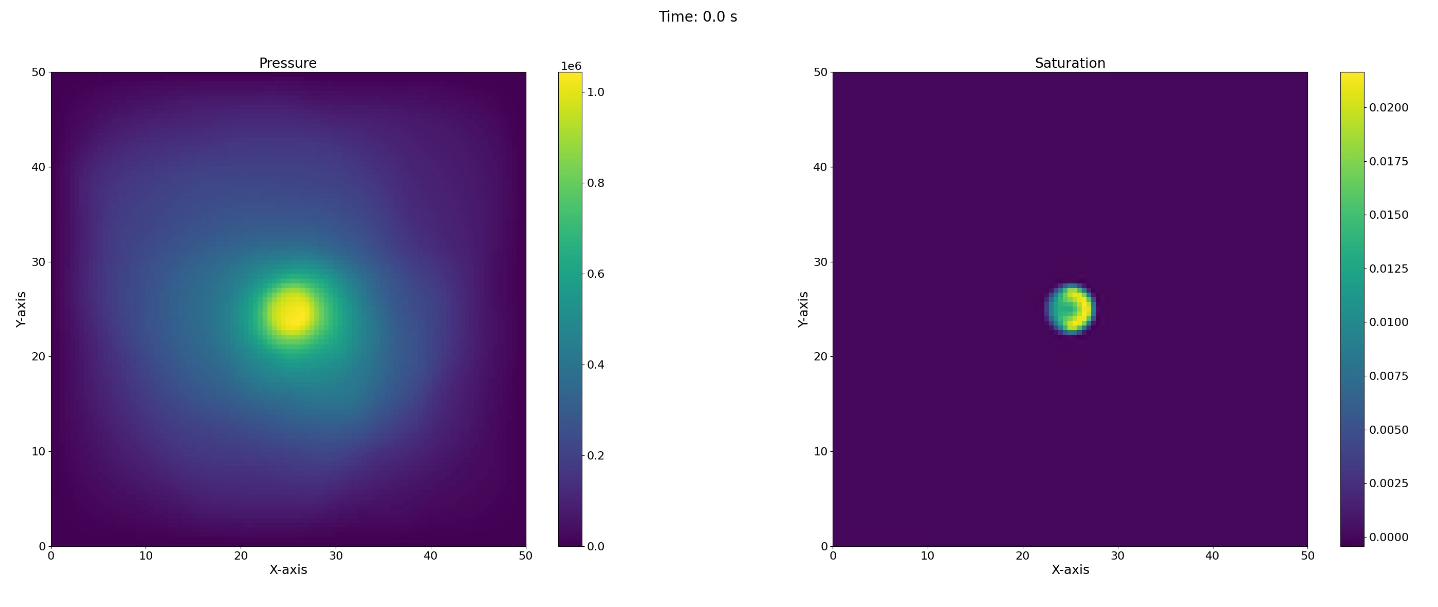}
        \caption{t = 0 s}
        \label{fig:K1_perm_map_t0}
    \end{subfigure}
    \hfill
    \begin{subfigure}[b]{0.9\textwidth}
        \centering
        \includegraphics[height=0.2\textheight]{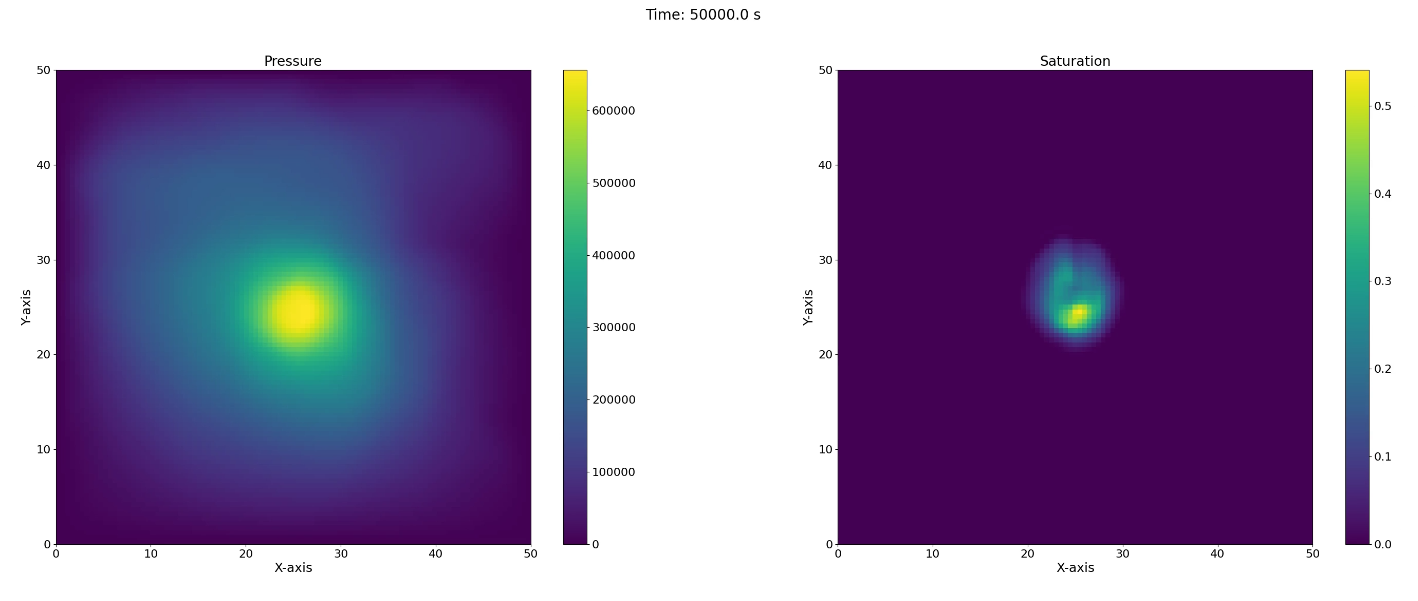}
        \caption{t = 50,000 s (= 50 time steps)}
        \label{fig:K1_perm_map_t100}
    \end{subfigure}
    \\
    \begin{subfigure}[b]{0.9\textwidth}
        \centering
        \includegraphics[height=0.2\textheight]{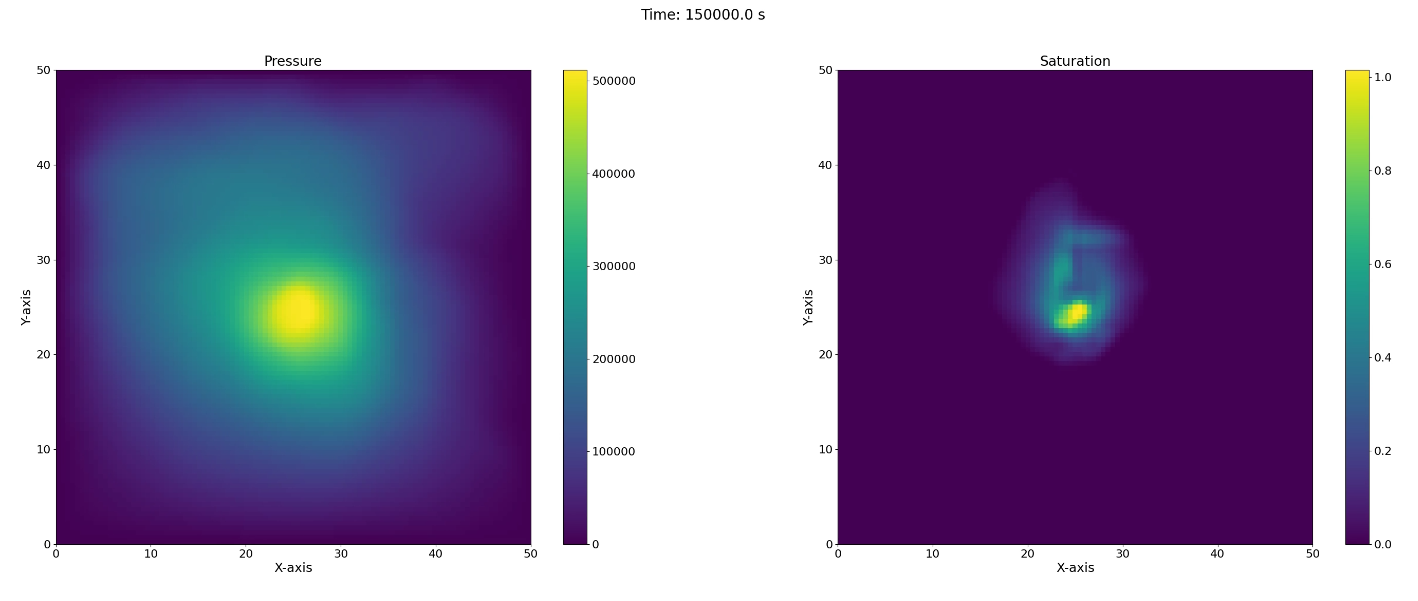}
        \caption{t = 150,000 s (= 150 time steps)}
        \label{fig:K1_perm_map_t200}
    \end{subfigure}
    \hfill
    \begin{subfigure}[b]{0.9\textwidth}
        \centering
        \includegraphics[height=0.2\textheight]{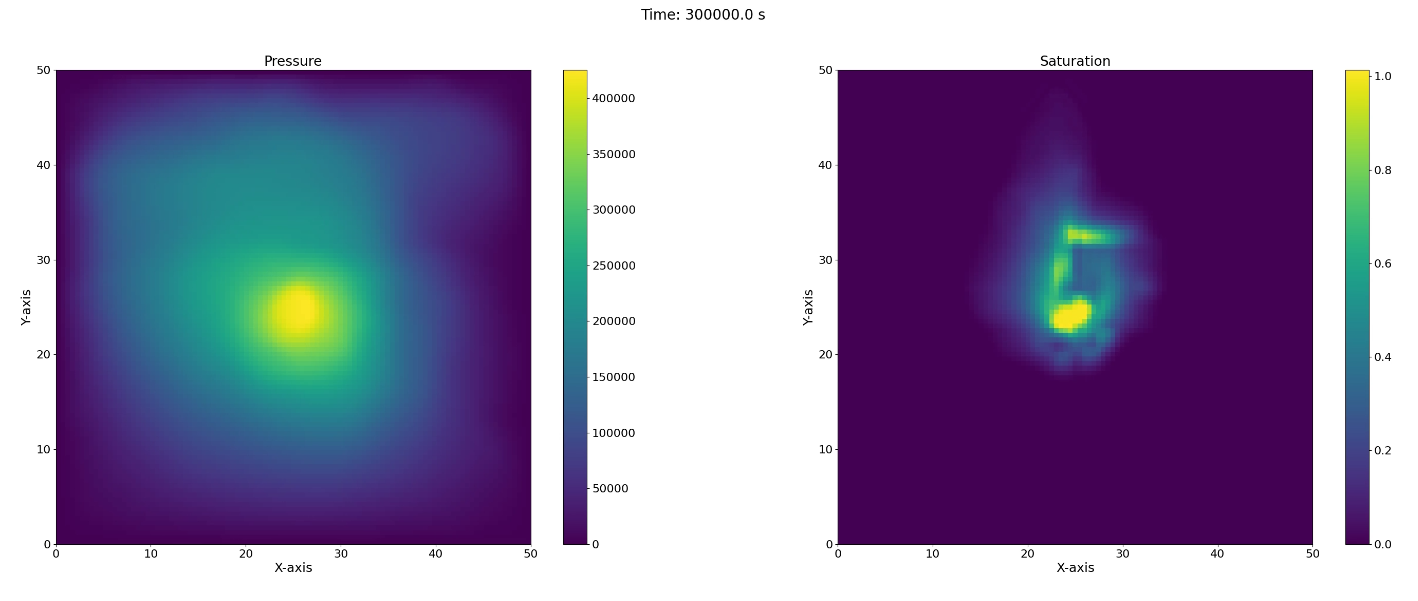}
        \caption{t = 300,000 s (= 300 time steps)} 
        \label{fig:K1_perm_map_t400}
    \end{subfigure}
    \caption{CO2 saturation evolution for configuration K1 with non-uniform permeability at different time steps.}
    \label{fig:results_perm_map_K1}
\end{figure}

\begin{figure}[htbp]
    \centering
    \begin{subfigure}[b]{0.9\textwidth}
        \centering
        \includegraphics[height=0.2\textheight]{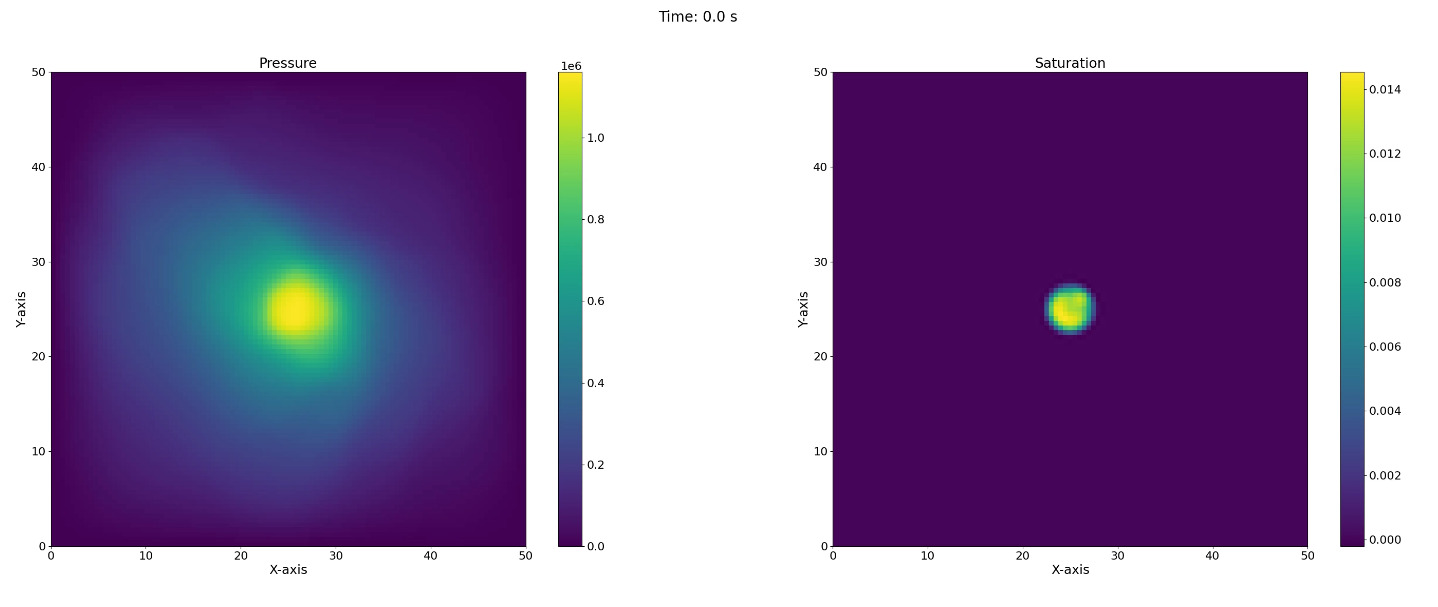}
        \caption{t = 0 s}
        \label{fig:K2_perm_map_t0}
    \end{subfigure}
    \hfill
    \begin{subfigure}[b]{0.9\textwidth}
        \centering
        \includegraphics[height=0.2\textheight]{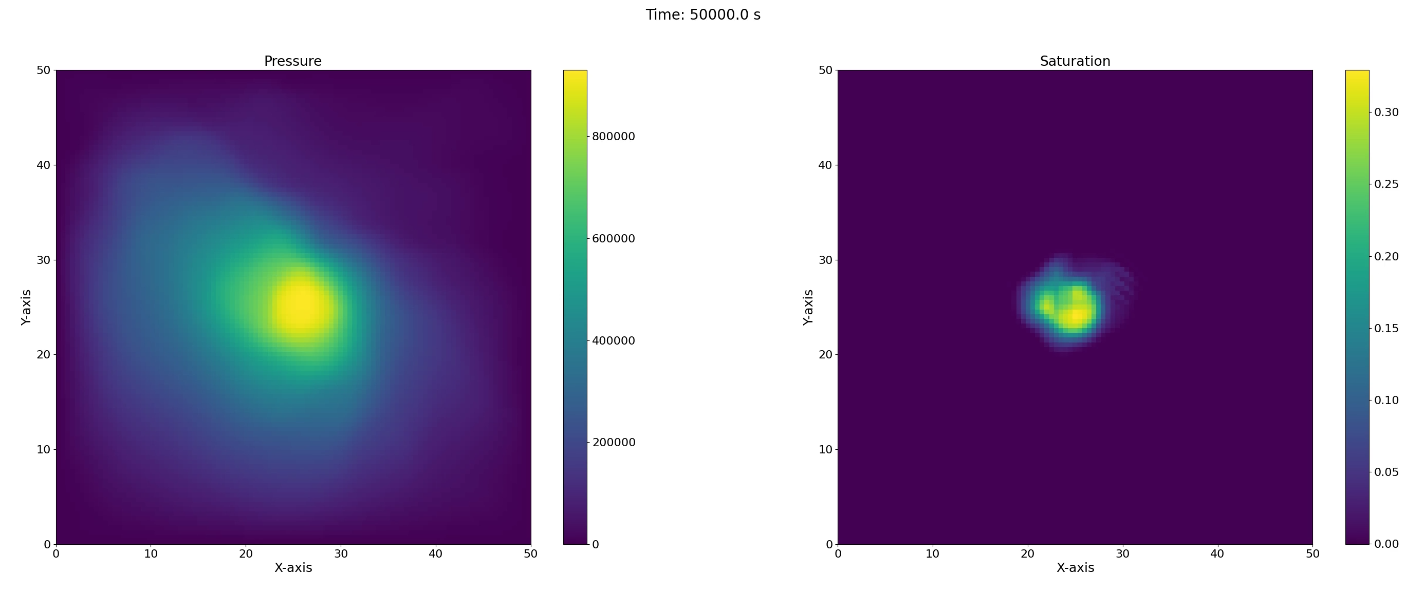}
        \caption{t = 50,000 s (= 50 time steps)}
        \label{fig:K2_perm_map_t100}
    \end{subfigure}
    \\
    \begin{subfigure}[b]{0.9\textwidth}
        \centering
        \includegraphics[height=0.2\textheight]{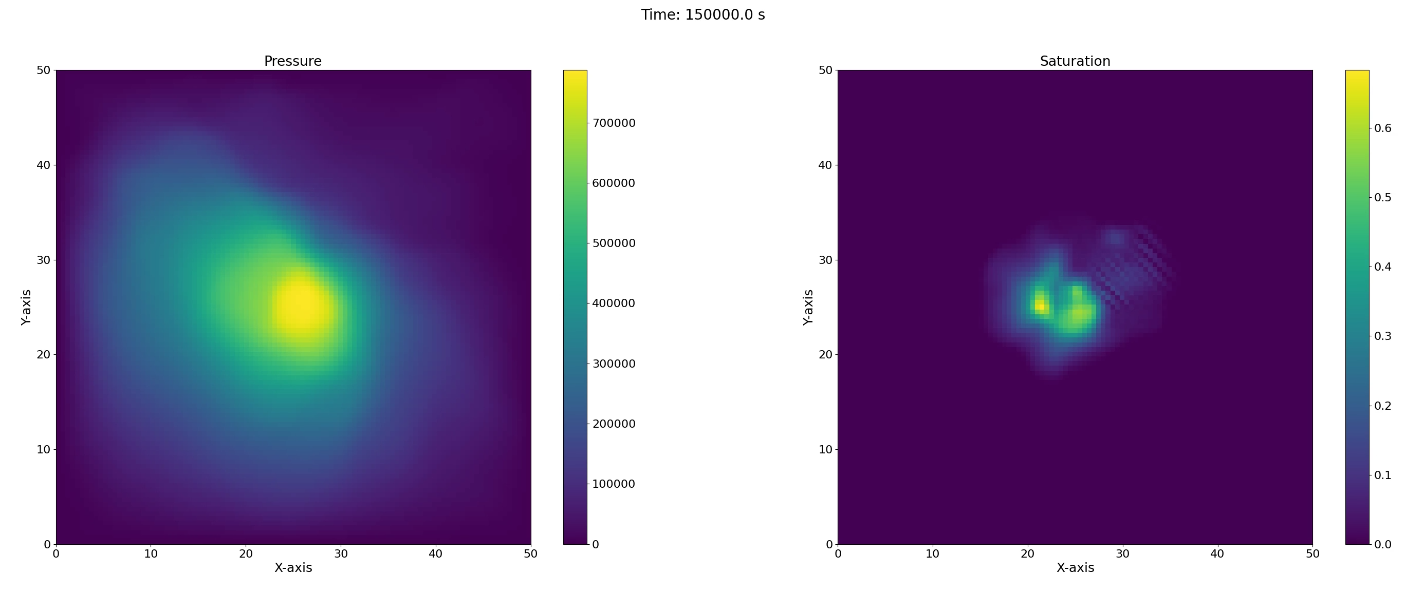}
        \caption{t = 150,000 s (= 150 time steps)}
        \label{fig:K2_perm_map_t200}
    \end{subfigure}
    \hfill
    \begin{subfigure}[b]{0.9\textwidth}
        \centering
        \includegraphics[height=0.2\textheight]{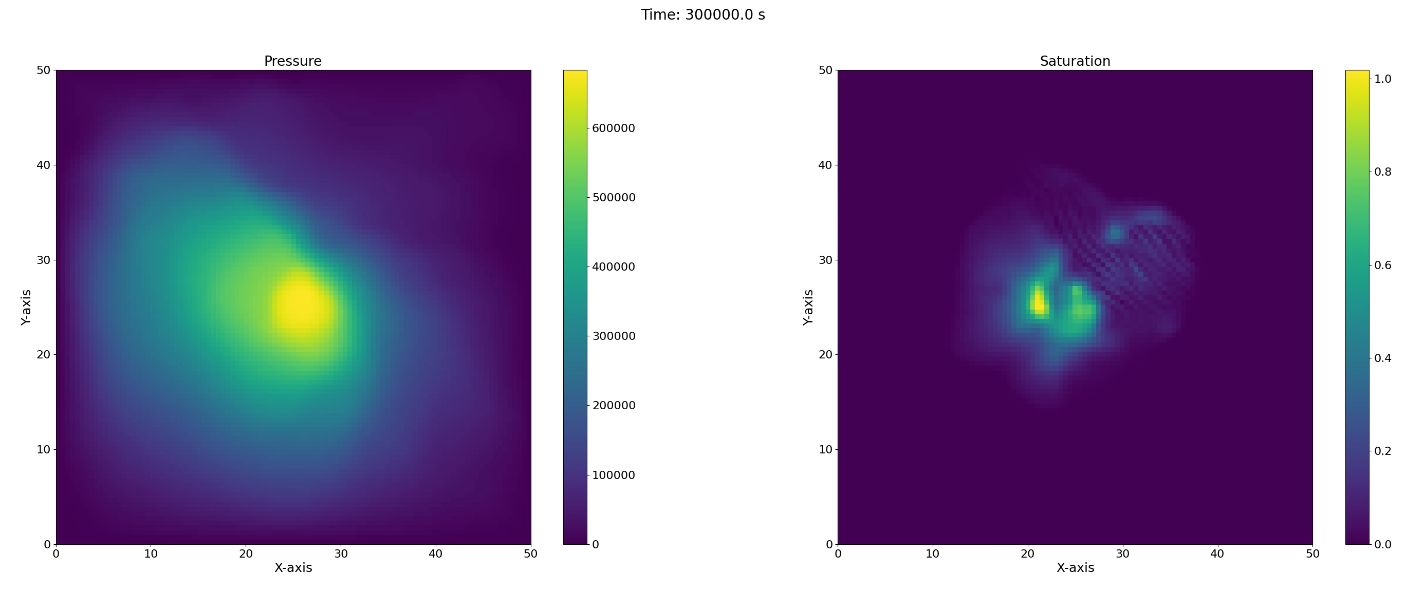}
        \caption{t = 300,000 s (= 300 time steps)}
        \label{fig:K2_perm_map_t400}
    \end{subfigure}
    \caption{CO2 saturation evolution for configuration K2 with non-uniform permeability at different time steps.}
    \label{fig:results_perm_map_K2}
\end{figure}

\begin{figure}[htbp]
    \centering
    \begin{subfigure}[b]{0.9\textwidth}
        \centering
        \includegraphics[height=0.2\textheight]{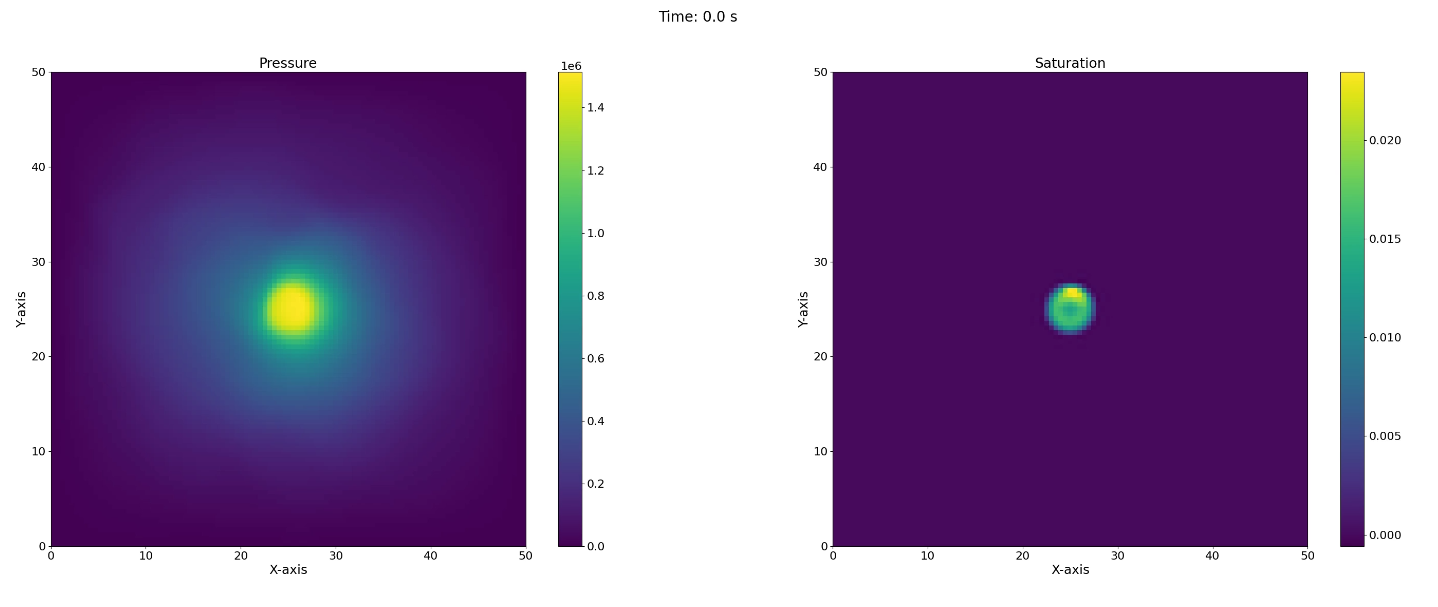}
        \caption{t = 0 s}
        \label{fig:K3_perm_map_t0}
    \end{subfigure}
    \hfill
    \begin{subfigure}[b]{0.9\textwidth}
        \centering
        \includegraphics[height=0.2\textheight]{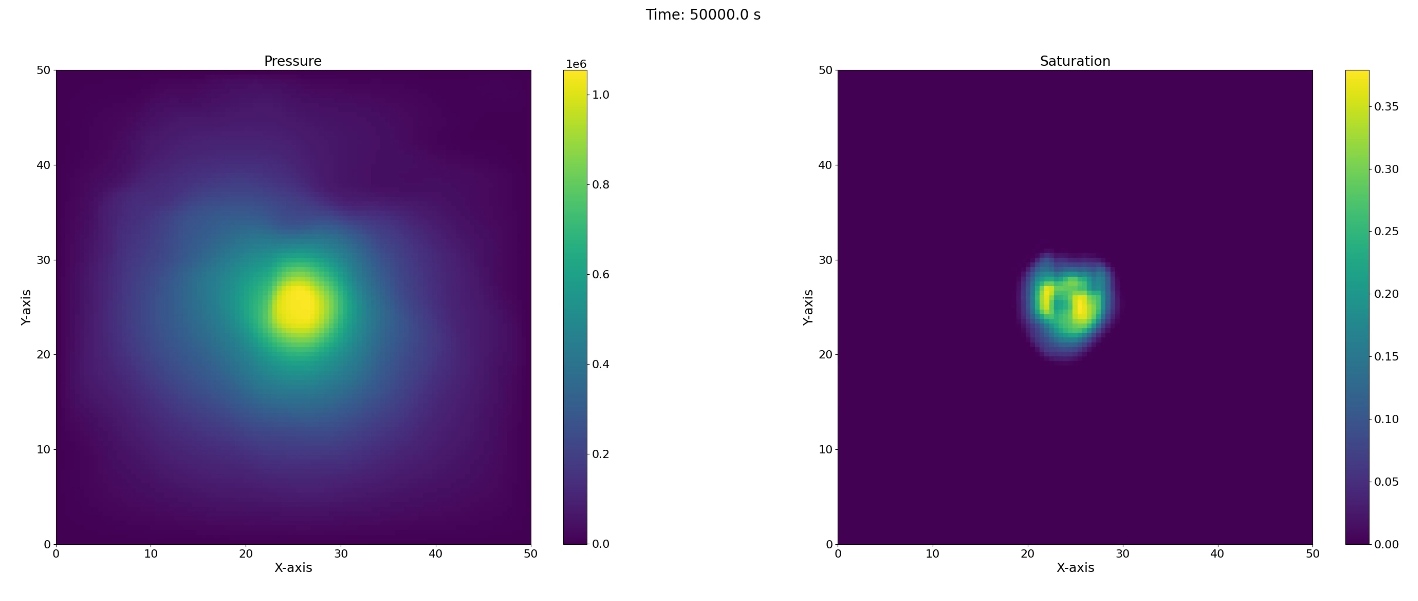}
        \caption{t = 50,000 s (= 50 time steps)}
        \label{fig:K3_perm_map_t100}
    \end{subfigure}
    \\
    \begin{subfigure}[b]{0.9\textwidth}
        \centering
        \includegraphics[height=0.2\textheight]{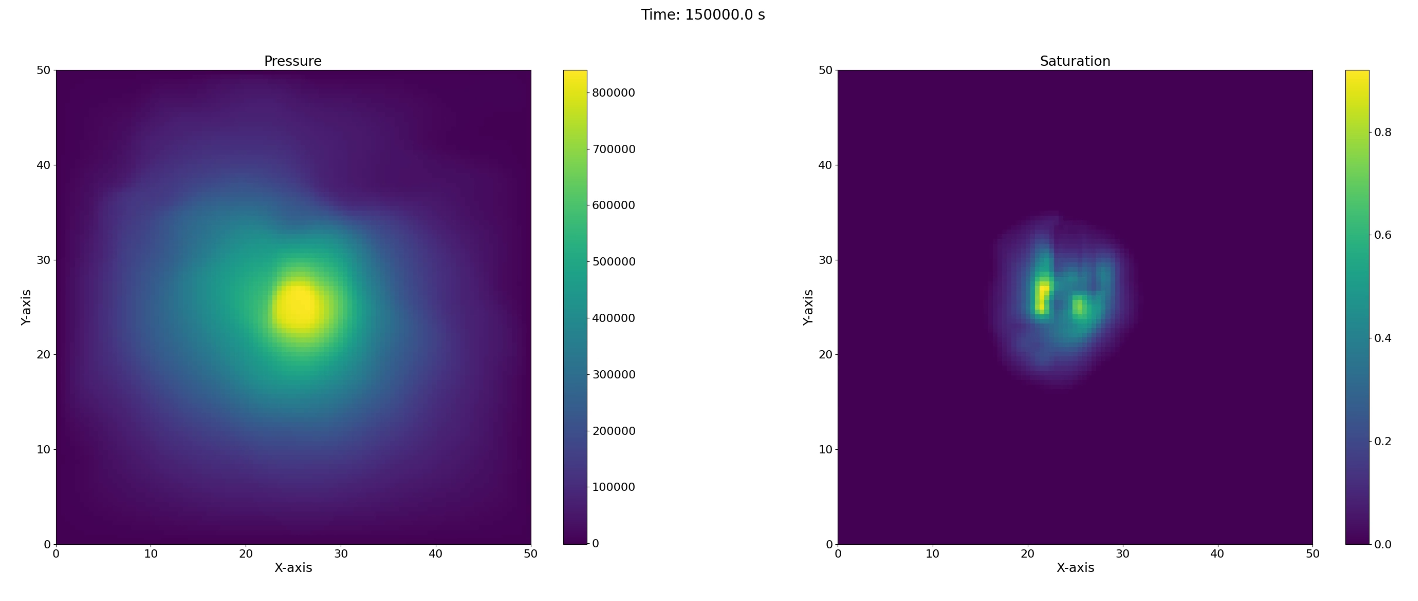}
        \caption{t = 150,000 s (= 150 time steps)}
        \label{fig:K3_perm_map_t200}
    \end{subfigure}
    \hfill
    \begin{subfigure}[b]{0.9\textwidth}
        \centering
        \includegraphics[height=0.2\textheight]{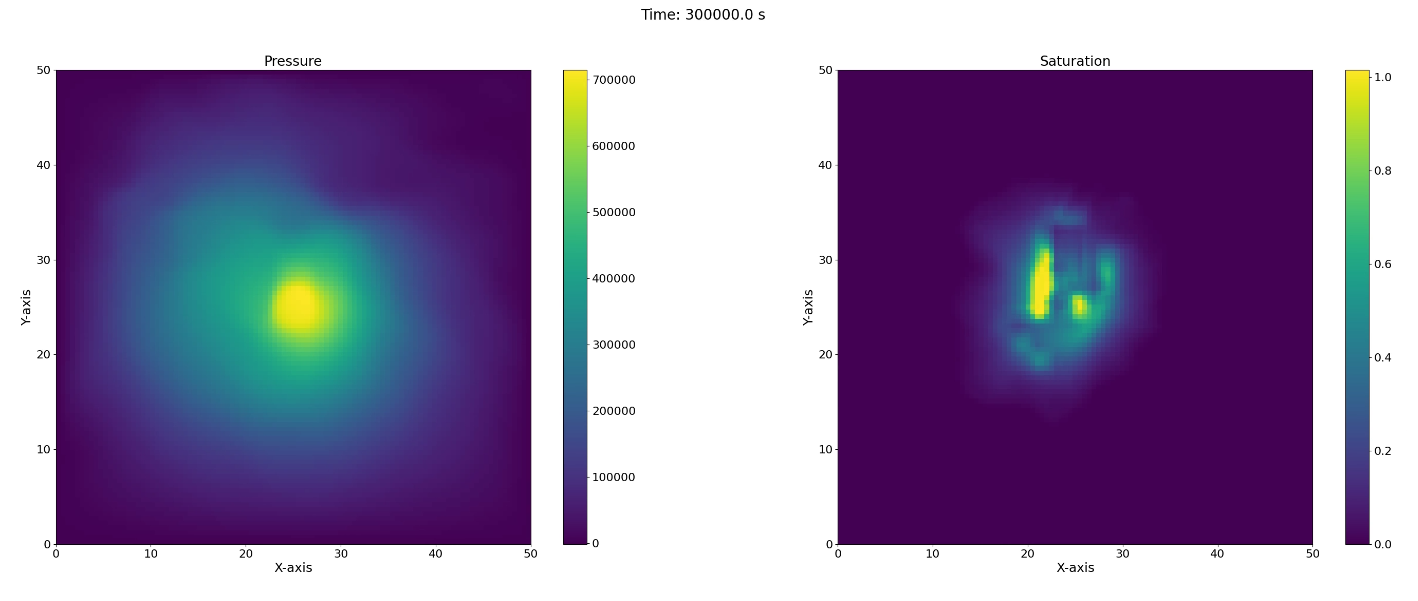}
        \caption{t = 300,000 s (= 300 time steps)}
        \label{fig:K3_perm_map_t400}
    \end{subfigure}
    \caption{CO2 saturation evolution for configuration K3 with non-uniform permeability at different time steps.}
    \label{fig:results_perm_map_K3}
\end{figure}

The next simulation set includes both porosity and permeability maps for the configurations K1, K2, and K3.
The simulation parameters are: time step size $\tau = 1000$ seconds, number of time steps = 500, and source term defined as a circle of radius 3 m at the center of the domain with a strength of $q_g = 5 \times 10^{-6}$ (dimensionless units of $S_g$) per iteration.
The results are shown in Figures \ref{fig:results_perm_map_K1}, \ref{fig:results_perm_map_K2}, \ref{fig:results_perm_map_K3}.

\subsection{Numerical experiments with IGA-ADS-CRVPINN solver}

\begin{figure}[htbp]
    \centering
    \begin{subfigure}[b]{0.45\textwidth}
        \centering
        \includegraphics[width=\textwidth]{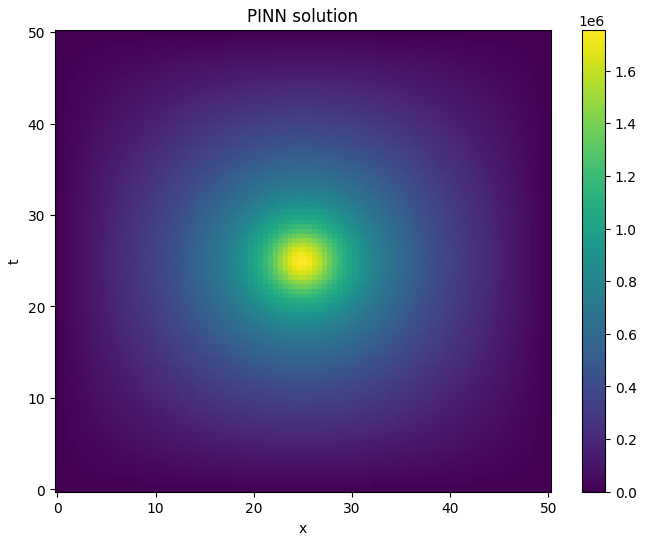}
        \caption{CRVPINN solution at t = 0 s}
        \label{fig:crvpinn_K1_t0}
    \end{subfigure}
    \hfill
    \begin{subfigure}[b]{0.45\textwidth}
        \centering
        \includegraphics[width=\textwidth]{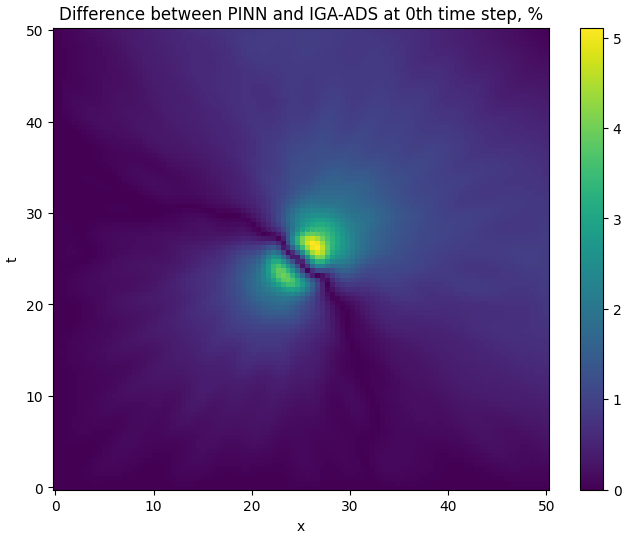}
        \caption{Difference at t = 0 s}
        \label{fig:crvpinn_K1_diff_t0}
    \end{subfigure}
    \\
    \begin{subfigure}[b]{0.45\textwidth}
        \centering
        \includegraphics[width=\textwidth]{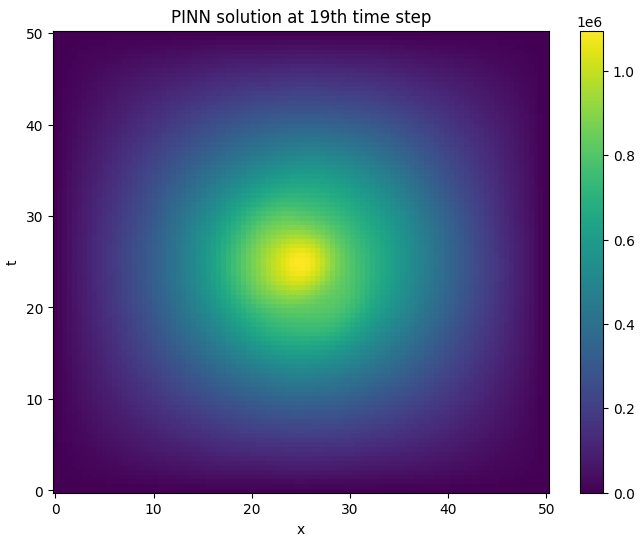}
        \caption{CRVPINN solution after 20 time steps}
        \label{fig:crvpinn_K1_t19}
    \end{subfigure}
    \hfill
    \begin{subfigure}[b]{0.45\textwidth}
        \centering
        \includegraphics[width=\textwidth]{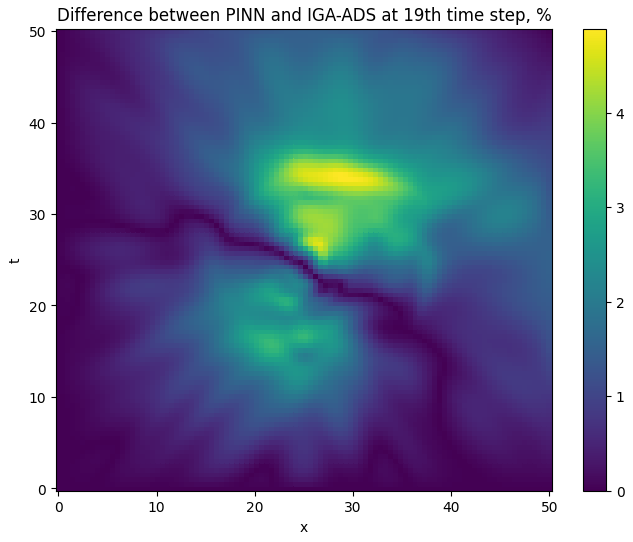}
        \caption{Difference after 20 time steps}
        \label{fig:crvpinn_K1_diff_t19}
    \end{subfigure}
    \caption{CRVPINN pressure solution comparison with IGA-ADS solver for configuration K1 with uniform permeability.}
    \label{fig:crvpinn_uniform_K1}
\end{figure}

\begin{figure}[htbp]
    \centering
    \begin{subfigure}[b]{0.45\textwidth}
        \centering
        \includegraphics[width=\textwidth]{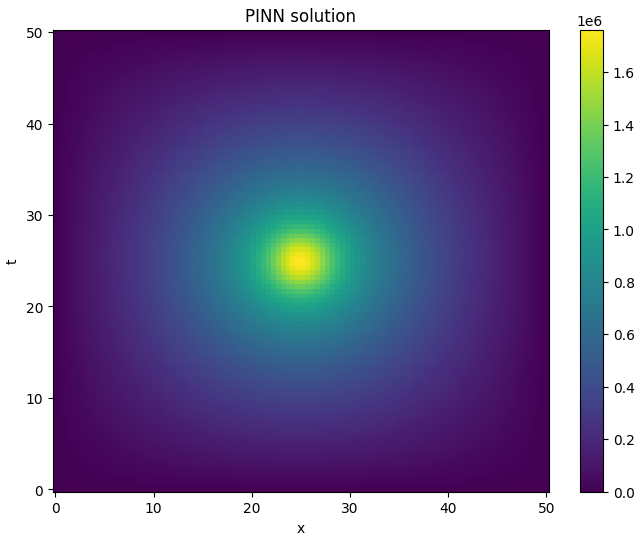}
        \caption{CRVPINN solution at t = 0 s}
        \label{fig:crvpinn_K2_t0}
    \end{subfigure}
    \hfill
    \begin{subfigure}[b]{0.45\textwidth}
        \centering
        \includegraphics[width=\textwidth]{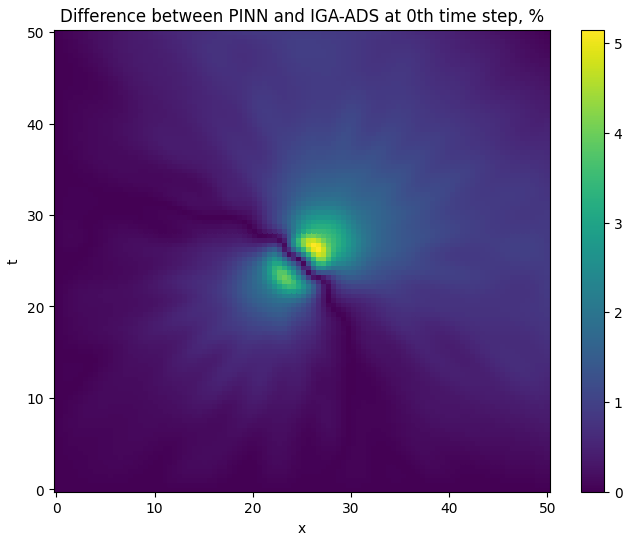}
        \caption{Difference at t = 0 s}
        \label{fig:crvpinn_K2_diff_t0}
    \end{subfigure}
    \\
    \begin{subfigure}[b]{0.45\textwidth}
        \centering
        \includegraphics[width=\textwidth]{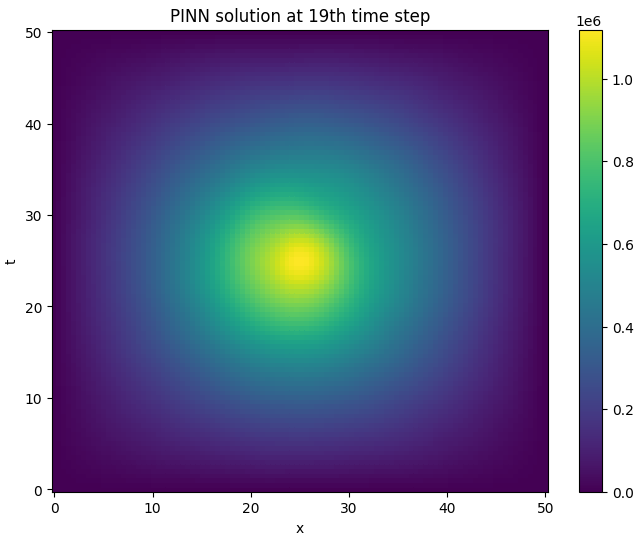}
        \caption{CRVPINN solution after 20 time steps}
        \label{fig:crvpinn_K2_t19}
    \end{subfigure}
    \hfill
    \begin{subfigure}[b]{0.45\textwidth}
        \centering
        \includegraphics[width=\textwidth]{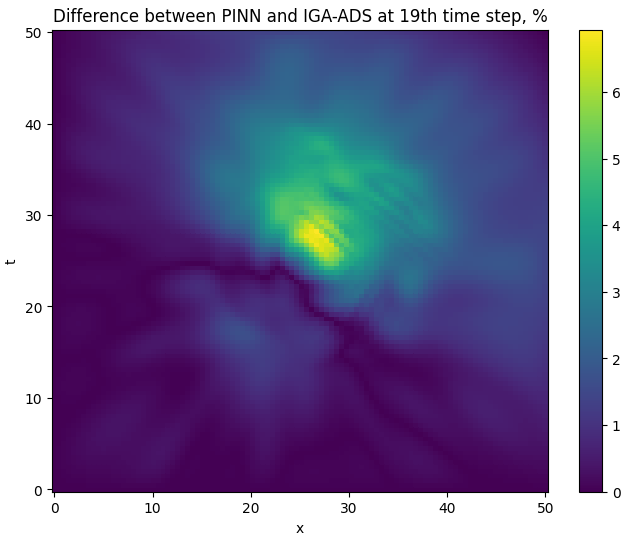}
        \caption{Difference after 20 time steps}
        \label{fig:crvpinn_K2_diff_t19}
    \end{subfigure}
    \caption{CRVPINN pressure solution comparison with IGA-ADS solver for configuration K2 with uniform permeability.}
    \label{fig:crvpinn_uniform_K2}
\end{figure}

\begin{figure}[htbp]
    \centering
    \begin{subfigure}[b]{0.45\textwidth}
        \centering
        \includegraphics[width=\textwidth]{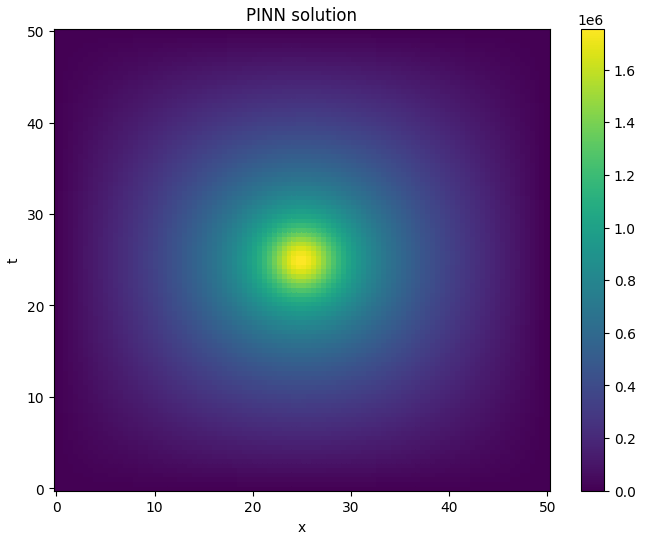}
        \caption{CRVPINN solution at t = 0 s}
        \label{fig:crvpinn_K3_t0}
    \end{subfigure}
    \hfill
    \begin{subfigure}[b]{0.45\textwidth}
        \centering
        \includegraphics[width=\textwidth]{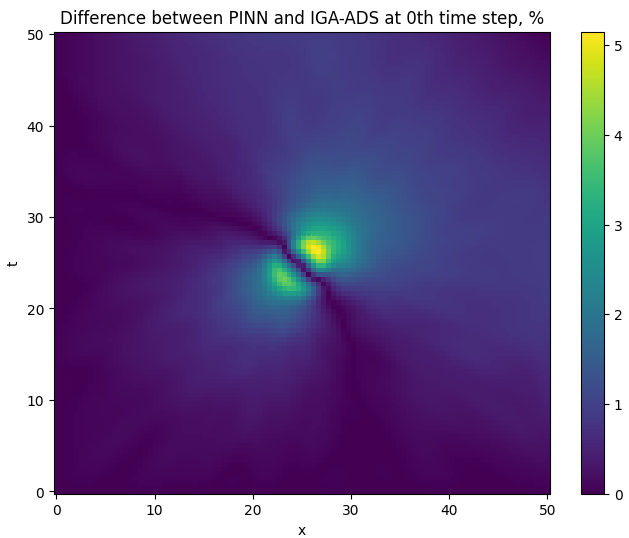}
        \caption{Difference at t = 0 s}
        \label{fig:crvpinn_K3_diff_t0}
    \end{subfigure}
    \\
    \begin{subfigure}[b]{0.45\textwidth}
        \centering
        \includegraphics[width=\textwidth]{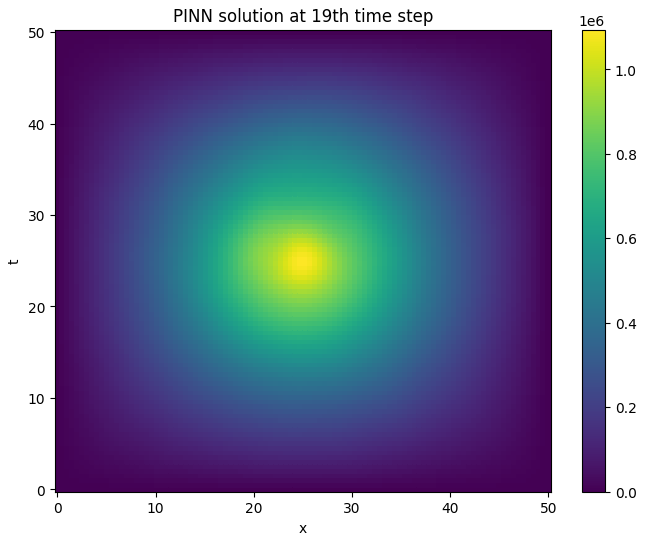}
        \caption{CRVPINN solution after 20 time steps}
        \label{fig:crvpinn_K3_t19}
    \end{subfigure}
    \hfill
    \begin{subfigure}[b]{0.45\textwidth}
        \centering
        \includegraphics[width=\textwidth]{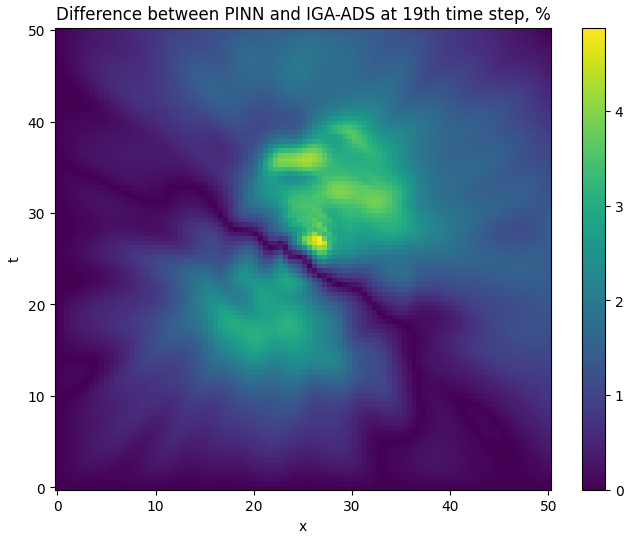}
        \caption{Difference after 20 time steps}
        \label{fig:crvpinn_K3_diff_t19}
    \end{subfigure}
    \caption{CRVPINN pressure solution comparison with IGA-ADS solver for configuration K3 with uniform permeability.}
    \label{fig:crvpinn_uniform_K3}
\end{figure}

Second, we performed the simulations of the IGA-ADS solver employed for the saturation updates \eqref{eq:6}, coupled with the CRVPINN solver equipped for the pressure updates \eqref{eq:5} (we use the IGA-ADS-CRVPINN solver).
We performed tests to validate the pressure evaluation against the IGA-ADS direct solver results based on the same simulation parameters. 
In our tests, we pretrained the CRVPINN model for the initial pressure configuration using 20,000 epochs with the Adam optimizer \cite{ADAM} and a learning rate of 0.0001 based on the initial saturation $S_g^0$.
In the figures, the initial training results are shown and compared with the IGA-ADS solver results.
After this pre-training step, we trained the model for each following time step for 100 epochs based on the new saturation $S_g^n$.
For the uniform permeability case, the results are shown in Figures \ref{fig:crvpinn_uniform_K1}, \ref{fig:crvpinn_uniform_K2}, \ref{fig:crvpinn_uniform_K3}.
We have been able to achieve very good agreement between the two solvers (with a relative error of less than 5\% in most of the domain).

\begin{figure}[htbp]
    \centering
    \begin{subfigure}[b]{0.45\textwidth}
        \centering
        \includegraphics[width=\textwidth]{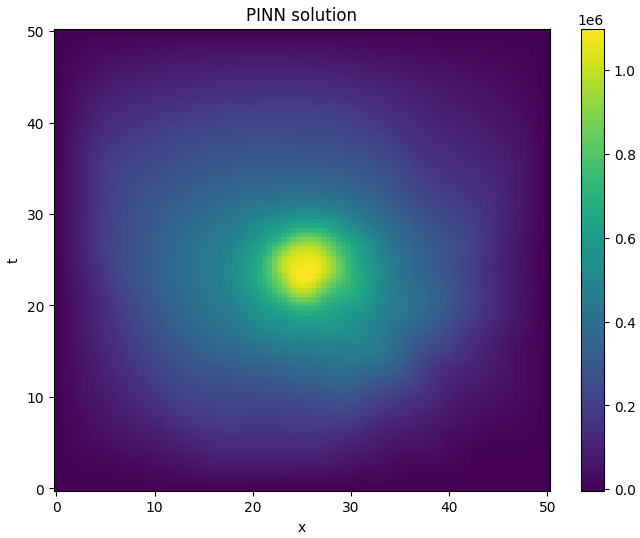}
        \caption{CRVPINN solution at t = 0 s}
        \label{fig:crvpinn_perm_K1_t0}
    \end{subfigure}
    \hfill
    \begin{subfigure}[b]{0.45\textwidth}
        \centering
        \includegraphics[width=\textwidth]{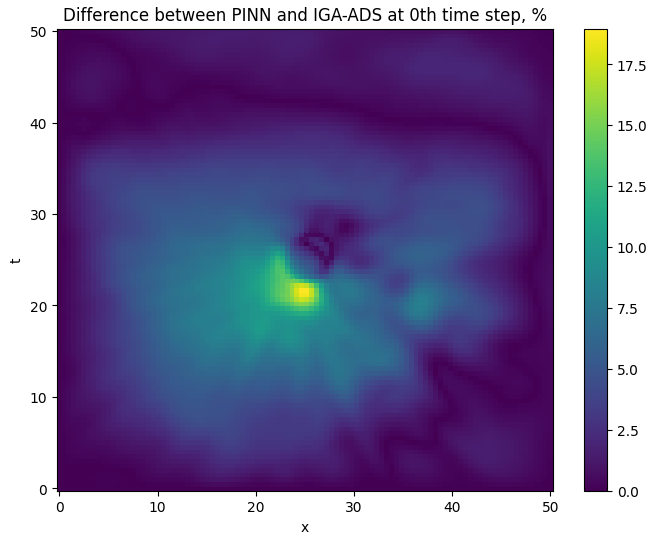}
        \caption{Difference at t = 0 s}
        \label{fig:crvpinn_perm_K1_diff_t0}
    \end{subfigure}
    \\
    \begin{subfigure}[b]{0.45\textwidth}
        \centering
        \includegraphics[width=\textwidth]{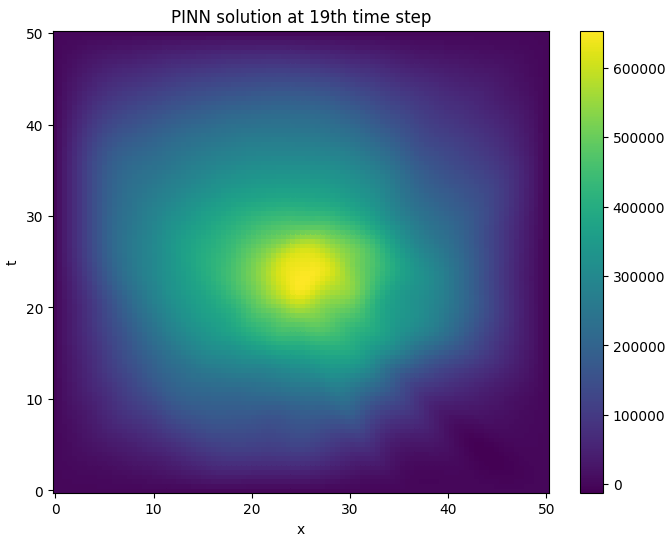}
        \caption{CRVPINN solution after 20 time steps}
        \label{fig:crvpinn_perm_K1_t19}
    \end{subfigure}
    \hfill
    \begin{subfigure}[b]{0.45\textwidth}
        \centering
        \includegraphics[width=\textwidth]{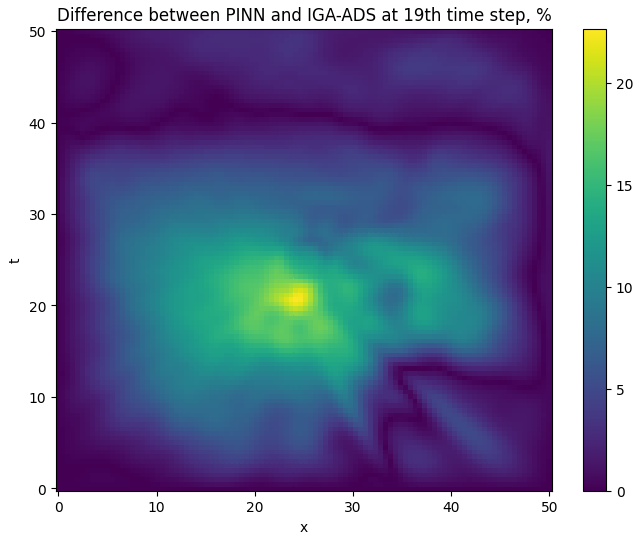}
        \caption{Difference after 20 time steps}
        \label{fig:crvpinn_perm_K1_diff_t19}
    \end{subfigure}
    \caption{CRVPINN pressure solution comparison with IGA-ADS solver for configuration K1 with permeability map.}
    \label{fig:crvpinn_perm_map_K1}
\end{figure}

\begin{figure}[htbp]
    \centering
    \begin{subfigure}[b]{0.45\textwidth}
        \centering
        \includegraphics[width=\textwidth]{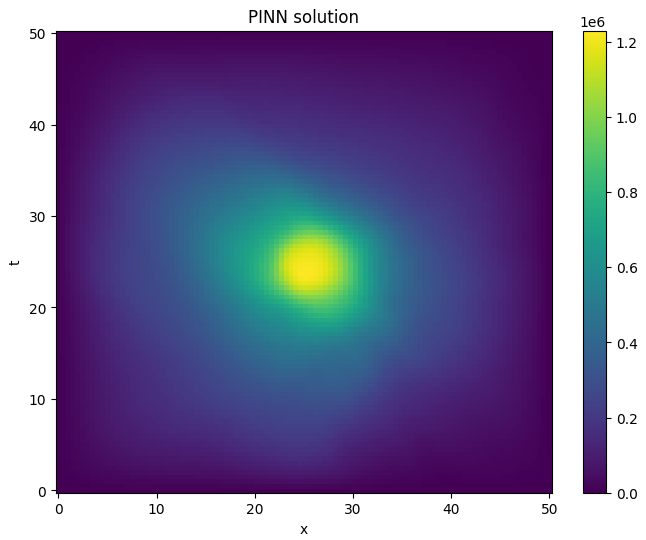}
        \caption{CRVPINN solution at t = 0 s}
        \label{fig:crvpinn_perm_K2_t0}
    \end{subfigure}
    \hfill
    \begin{subfigure}[b]{0.45\textwidth}
        \centering
        \includegraphics[width=\textwidth]{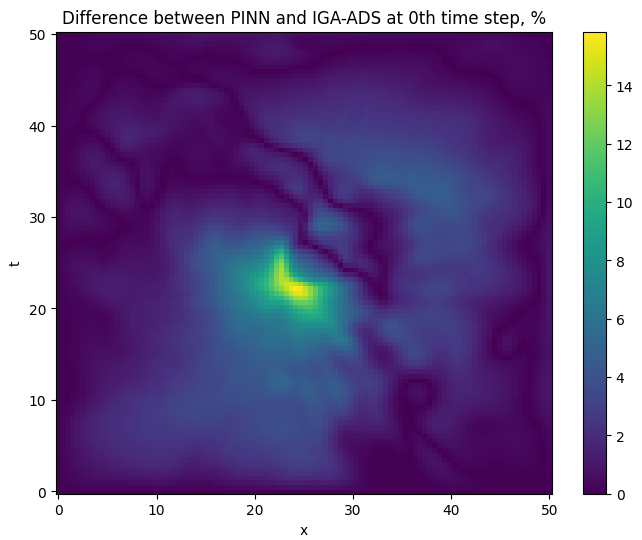}
        \caption{Difference at t = 0 s}
        \label{fig:crvpinn_perm_K2_diff_t0}
    \end{subfigure}
    \\
    \begin{subfigure}[b]{0.45\textwidth}
        \centering
        \includegraphics[width=\textwidth]{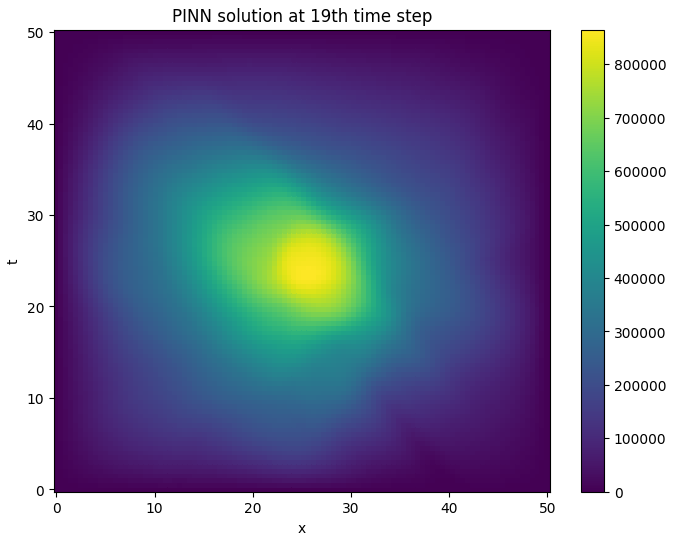}
        \caption{CRVPINN solution after 20 time steps}
        \label{fig:crvpinn_perm_K2_t19}
    \end{subfigure}
    \hfill
    \begin{subfigure}[b]{0.45\textwidth}
        \centering
        \includegraphics[width=\textwidth]{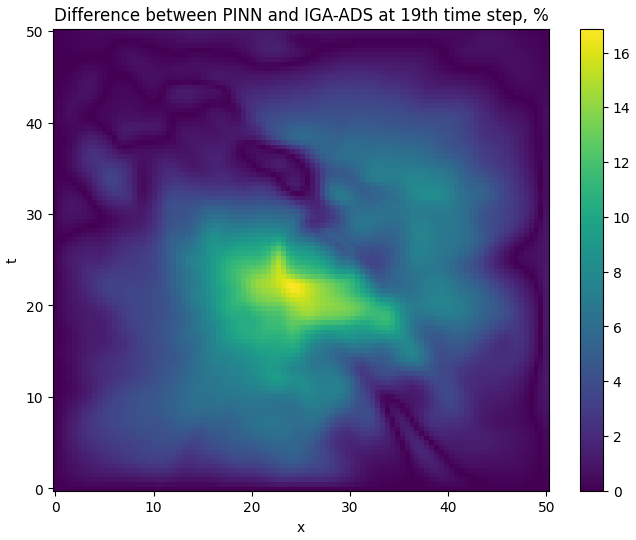}
        \caption{Difference after 20 time steps}
        \label{fig:crvpinn_perm_K2_diff_t19}
    \end{subfigure}
    \caption{CRVPINN pressure solution comparison with IGA-ADS solver for configuration K2 with permeability map.}
    \label{fig:crvpinn_perm_map_K2}
\end{figure}

\begin{figure}[htbp]
    \centering
    \begin{subfigure}[b]{0.45\textwidth}
        \centering
        \includegraphics[width=\textwidth]{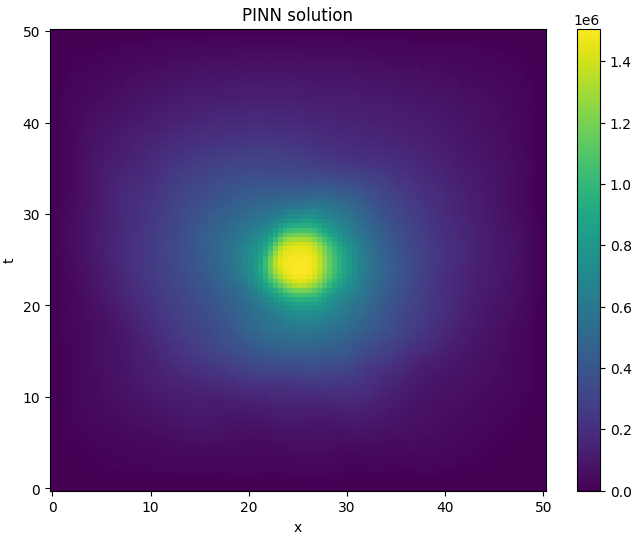}
        \caption{CRVPINN solution at t = 0 s}
        \label{fig:crvpinn_perm_K3_t0}
    \end{subfigure}
    \hfill
    \begin{subfigure}[b]{0.45\textwidth}
        \centering
        \includegraphics[width=\textwidth]{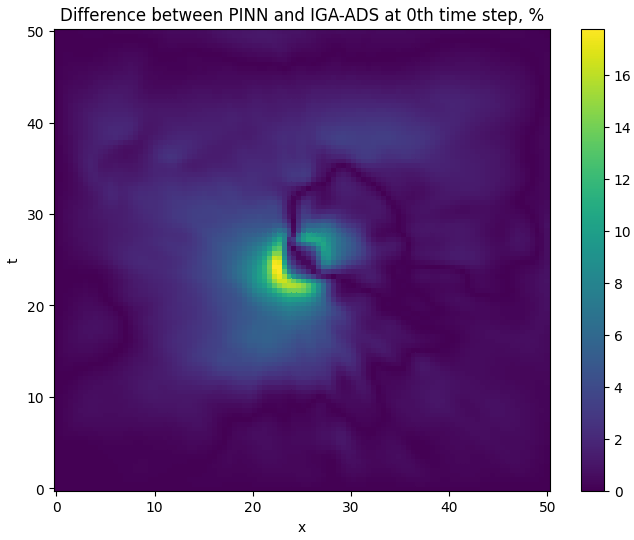}
        \caption{Difference at t = 0 s}
        \label{fig:crvpinn_perm_K3_diff_t0}
    \end{subfigure}
    \\
    \begin{subfigure}[b]{0.45\textwidth}
        \centering
        \includegraphics[width=\textwidth]{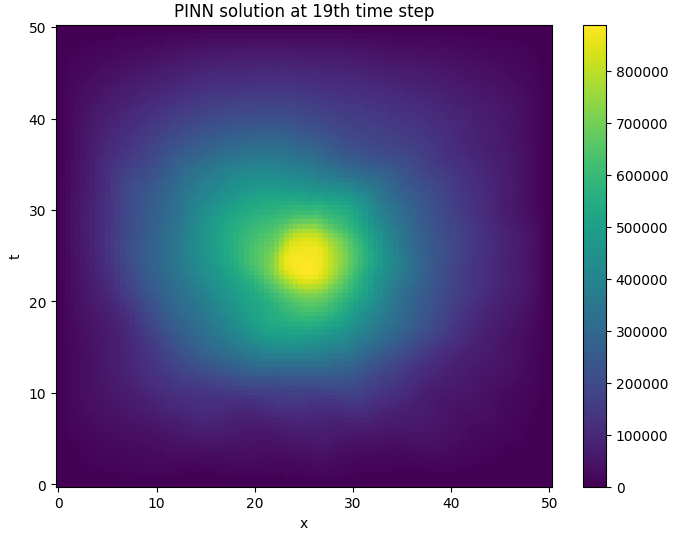}
        \caption{CRVPINN solution after 20 time steps}
        \label{fig:crvpinn_perm_K3_t19}
    \end{subfigure}
    \hfill
    \begin{subfigure}[b]{0.45\textwidth}
        \centering
        \includegraphics[width=\textwidth]{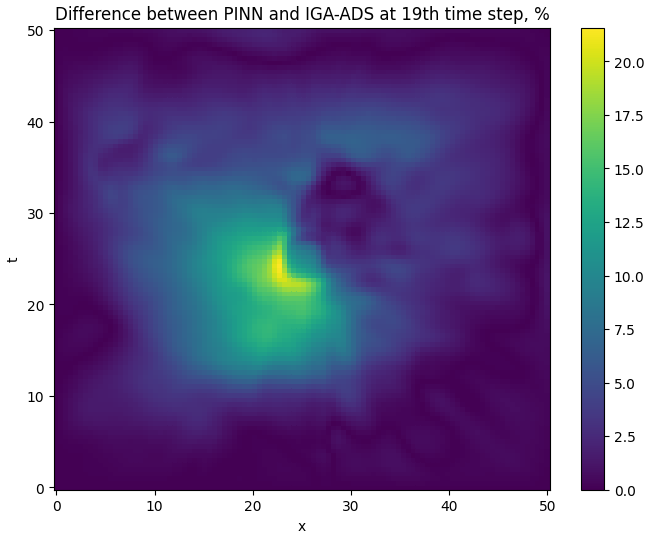}
        \caption{Difference after 20 time steps}
        \label{fig:crvpinn_perm_K3_diff_t19}
    \end{subfigure}
    \caption{CRVPINN pressure solution comparison with IGA-ADS solver for configuration K3 with permeability map.}
    \label{fig:crvpinn_perm_map_K3}
\end{figure}

For non-uniform permeability cases, the results are shown in Figures \ref{fig:crvpinn_perm_map_K1}, \ref{fig:crvpinn_perm_map_K2}, \ref{fig:crvpinn_perm_map_K3}.
The complexity of the permeability map introduces more challenges for the CRVPINN model, but the agreement with the IGA-ADS solver remains reasonable (with a relative error of about 15-20\% in some areas).
We have noticed that the model behaves better when the source term is larger in area or stronger in intensity.

\begin{figure}[ht]
    \centering
    \begin{minipage}{0.49\textwidth}
        \centering
        \begin{subfigure}[b]{0.9\textwidth}
            \centering
            \includegraphics[height=0.25\textheight]{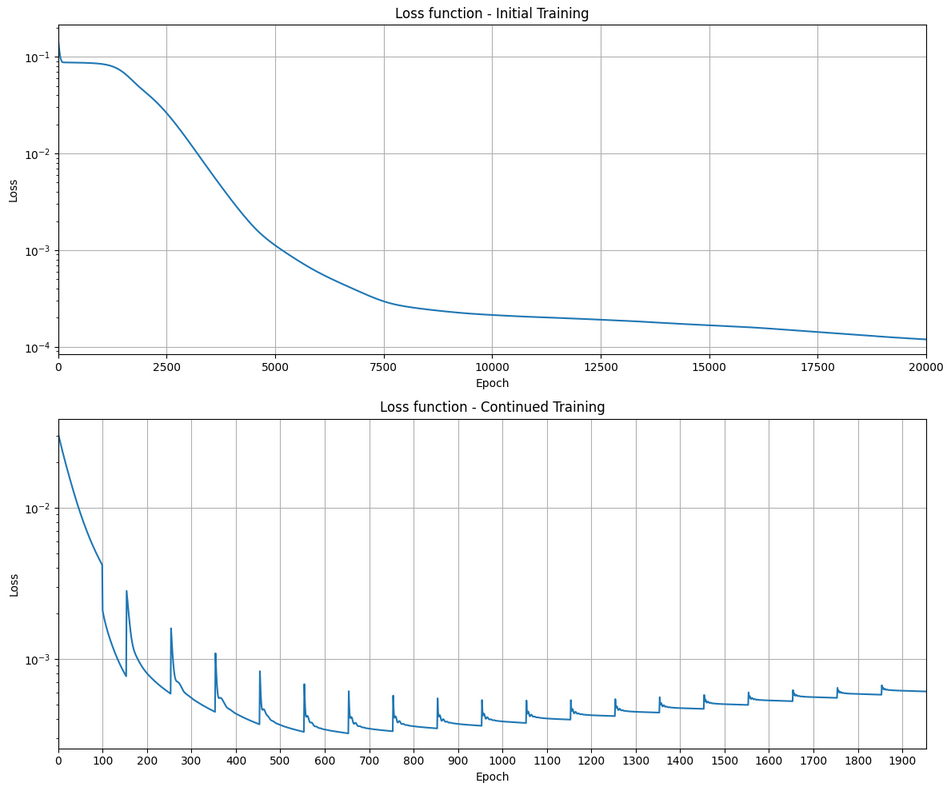}
            \caption{Configuration K1}
            \label{fig:loss_uniform_K1}
        \end{subfigure}
        \\
        \begin{subfigure}[b]{0.9\textwidth}
            \centering
            \includegraphics[height=0.25\textheight]{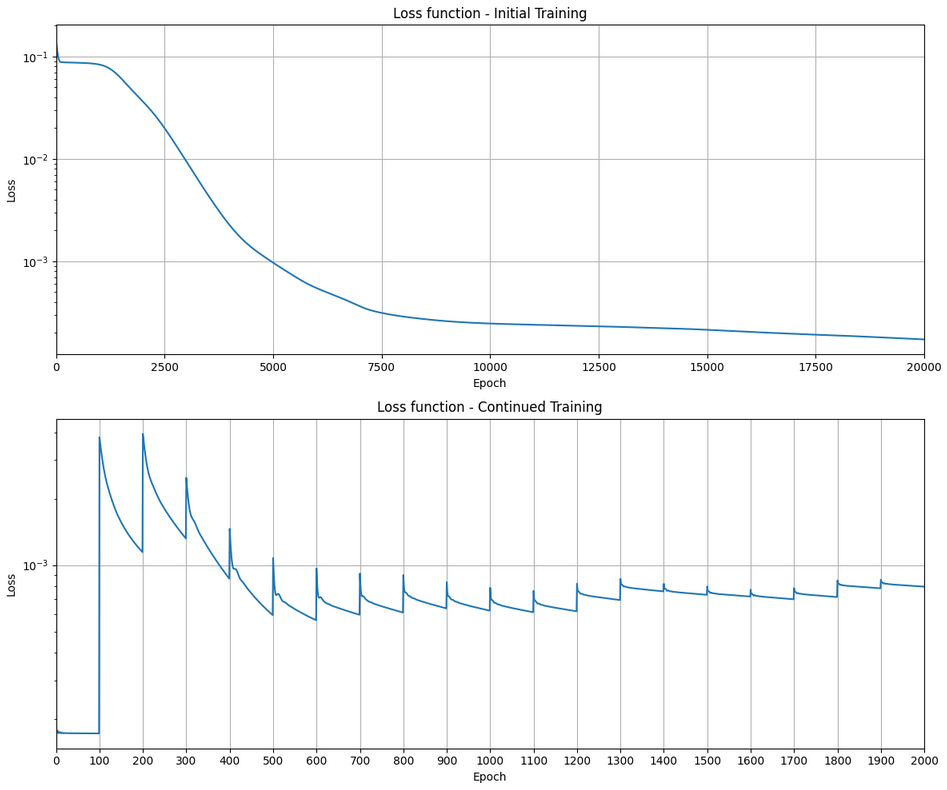}
            \caption{Configuration K2}
            \label{fig:loss_uniform_K2}
        \end{subfigure}
        \\
        \begin{subfigure}[b]{0.9\textwidth}
            \centering
            \includegraphics[height=0.25\textheight]{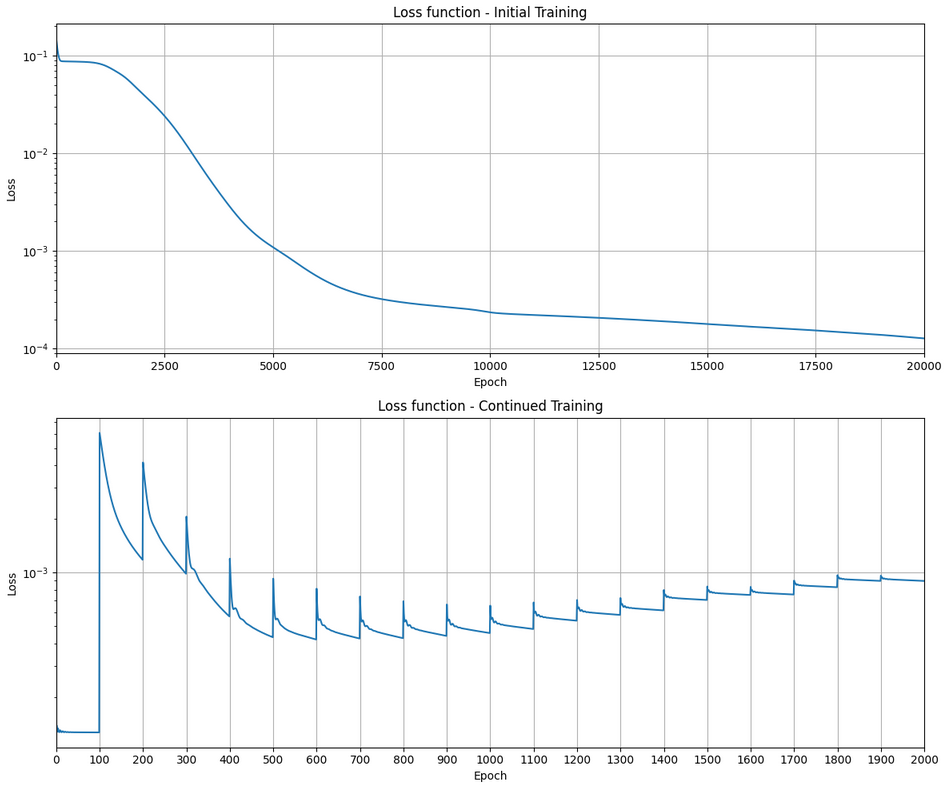}
            \caption{Configuration K3}
            \label{fig:loss_uniform_K3}
        \end{subfigure}
        \captionof{figure}{CRVPINN loss for uniform permeability field}
        \label{fig:uniform_loss}
    \end{minipage}
    \hfill
    \begin{minipage}{0.49\textwidth}
        \centering
        \begin{subfigure}[b]{0.9\textwidth}
            \centering
            \includegraphics[height=0.25\textheight]{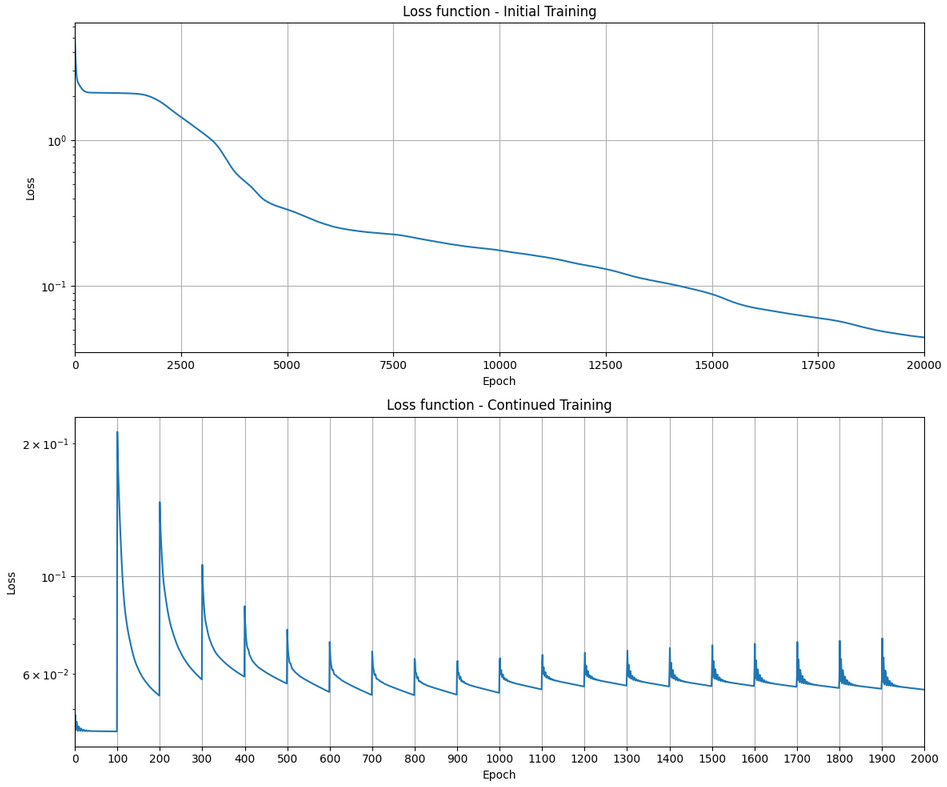}
            \caption{Configuration K1}
            \label{fig:loss_perm_map_K1}
        \end{subfigure}
        \\
        \begin{subfigure}[b]{0.9\textwidth}
            \centering
            \includegraphics[height=0.25\textheight]{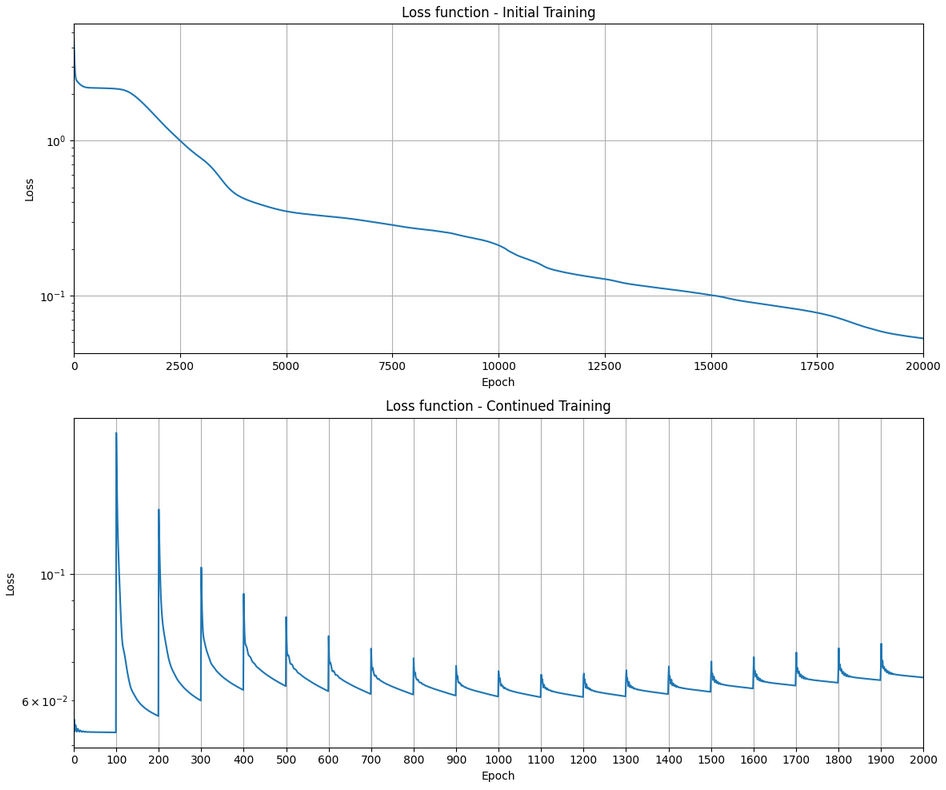}
            \caption{Configuration K2}
            \label{fig:loss_perm_map_K2}
        \end{subfigure}
        \\
        \begin{subfigure}[b]{0.9\textwidth}
            \centering
            \includegraphics[height=0.25\textheight]{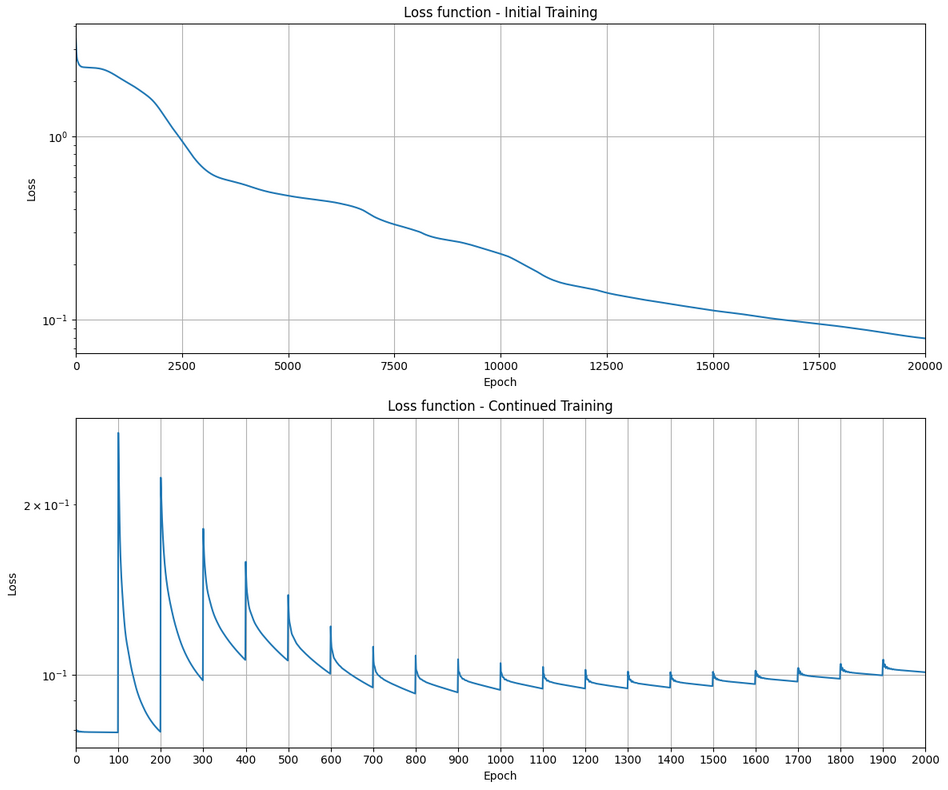}
            \caption{Configuration K3}
            \label{fig:loss_perm_map_K3}
        \end{subfigure}
        \captionof{figure}{CRVPINN loss for non-uniform permeability maps}
        \label{fig:perm_map_loss}
    \end{minipage}
\end{figure}

We have also recorded the loss function values during the training process for both uniform and non-uniform permeability cases. 
The loss curves are presented in Figures \ref{fig:uniform_loss} and \ref{fig:perm_map_loss}, respectively.
The top panel of each subfigure represents the 20,000 iterations of pretraining, and the bottom panel presents the pressure updates using 100 iterations of the training, for the first 200 iterations of the IGA-ADS (the CRVPINN training is executed every 10 iterations of the IGA-ADS solver).

\subsection{Timing and Performance Analysis}

For the training of the CRVPINN model, we use the Adam optimizer with a learning rate of $10^{-4}$, and a domain with $10^4$ collocation points (i.e. $N=100$), with the physical size of the grid being $50 \times 50$ m$^2$.
There are two main phases of the training process: pre-training and continuous training.
The pre-training phase is performed only once, and it consists of training the model for 20000 epochs on some starting configuration of the problem, given the corresponding saturation map from the IGA-ADS solver.
This phase is computationally expensive and takes around 500 to 600 seconds on a single GPGPU.
However, after this phase, we can perform continuous training for the subsequent time steps, which takes around 3-5 seconds for each step (100 epochs every 10 iterations from the IGA-ADS solver).
For each of these time steps, we use the saturation map from the IGA-ADS solver as input for the CRVPINN model and then return a pressure map, which is used for the next iteration of the IGA-ADS solver.

\begin{figure}
\includegraphics[width=\textwidth]{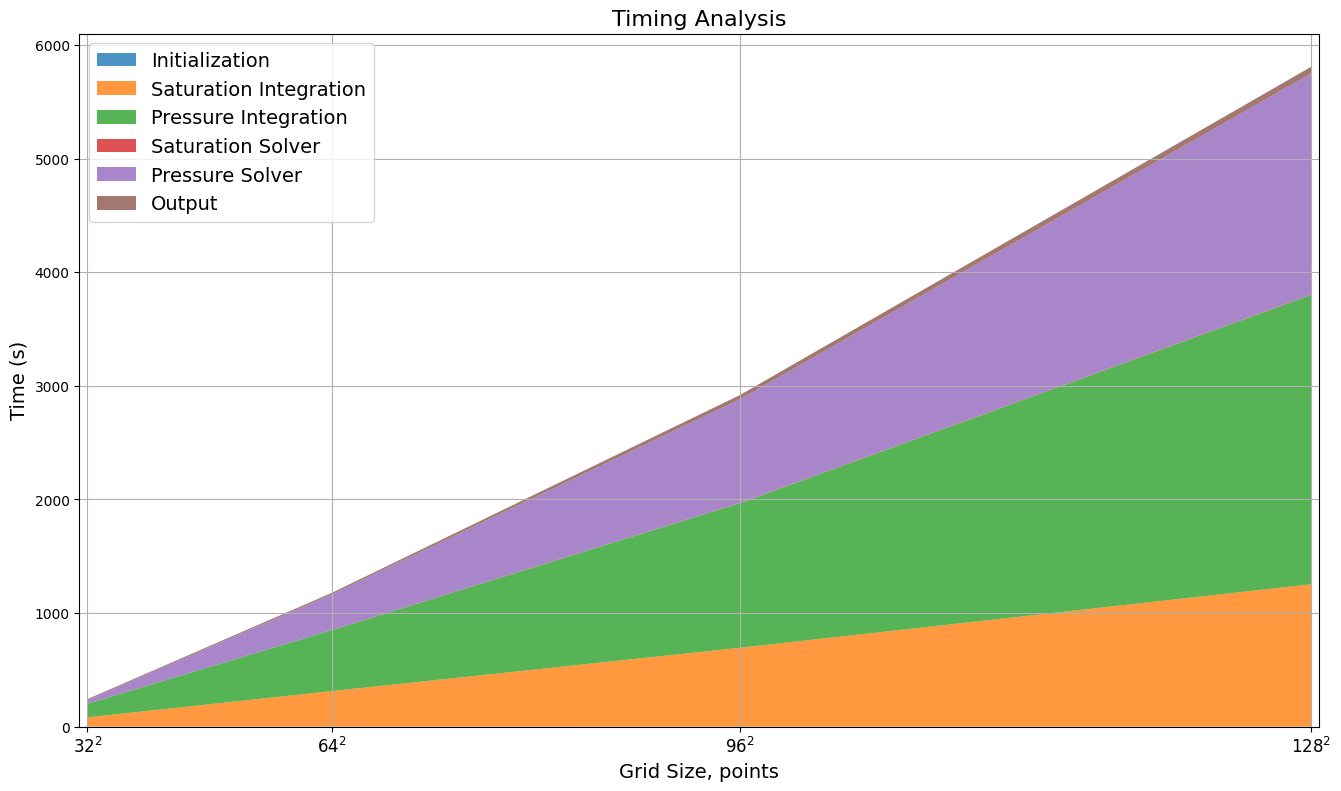}
\caption{Execution times for IGA-ADS solver for saturation coupled with MUMPS solver executed every 10 time steps for the pressure. Mesh dimensions for both IGA-ADS and CRVPINN vary from $32^2$, $64^2$, $96^2$, up to $128^2$. Quadratic B-splines with $C^1$ continuity. }
\label{fig:IGAMUMPS}
\end{figure}

\begin{figure}
\includegraphics[width=\textwidth]{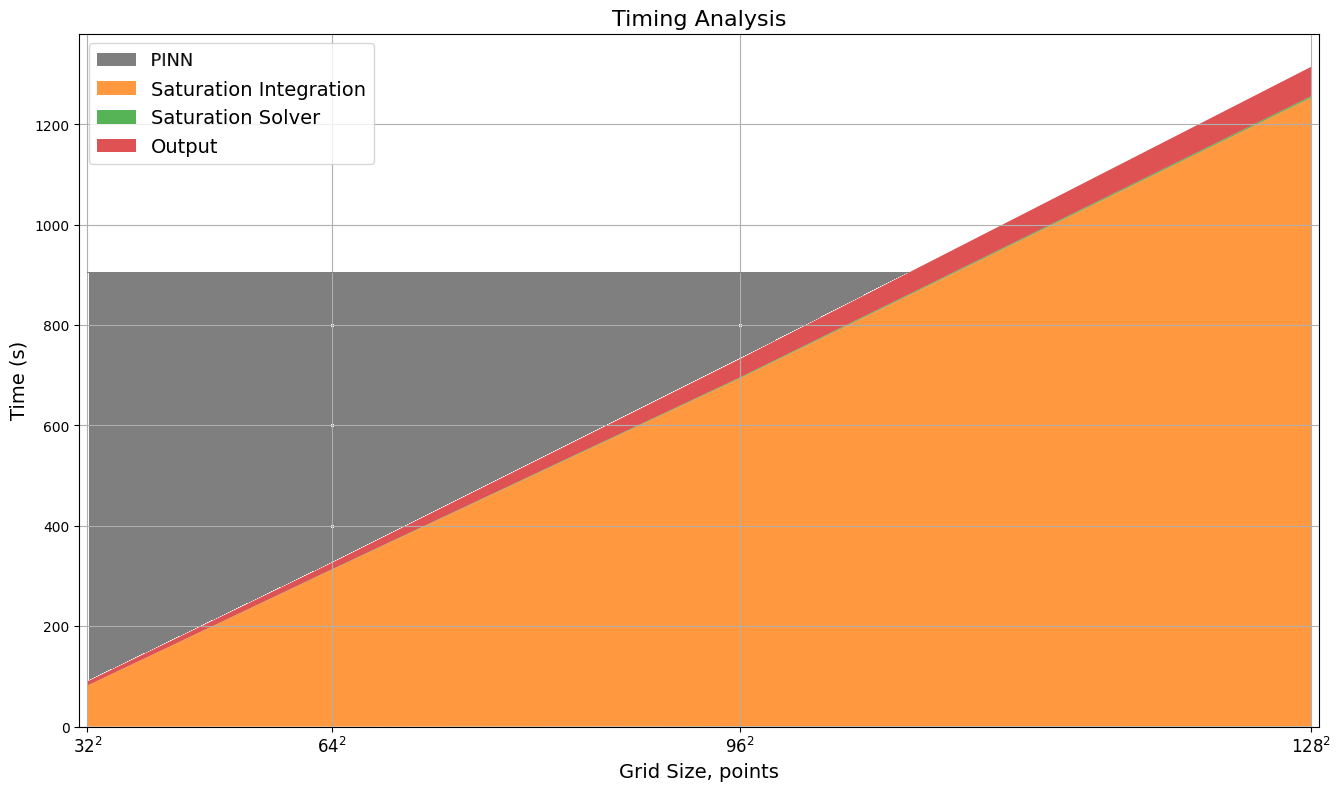}
\caption{Execution times for IGA-ADS solver for saturation coupled with CVPINN solver executed every 10 time steps for the pressure. Mesh dimensions for IGA-ADS varying from $32^2$, $64^2$, $96^2$, up to $128^2$. Quadratic B-splines with $C^1$ continuity. The number of collocation points for the CRVPINN solver is always $100^2$. }
\label{fig:IGACRVPINN}
\end{figure}

A comparison between the IGA-ADS solver for saturation, coupled with the MUMPS direct solver for pressure, and the IGA-ADS solver for saturation, coupled with the CRVPINN solver for pressure, is presented in Figures \ref{fig:IGAMUMPS}-\ref{fig:IGACRVPINN}. 
It has been developed on a computing node from the ARES cluster \cite{Ares} at ACC CYFRONET\footnote{Academic {C}omputer {C}entre {C}yfronet {A}{G}{H}  https://www.cyfronet.pl
}. The node is equipped with 
32 cores, an Intel(R) Xeon(R) Gold 6242 CPU @ 2.80GHz, and 384 GB of RAM.
% Askold: GPU information is irrelevant for ARES, we only use it for CPU computations
% all the GPU computations (CRVPINN) were computed on ATHENA.
% along with an NVIDIA Tesla V100-SXM2 GPGPU.
We executed 1000 time steps with the IGA-ADS-MUMPS solver. The timings are presented in Figure \ref{fig:IGAMUMPS}. Both the IGA-ADS solver and the MUMPS direct solver were executed on a single workstation of the CYFRONET cluster using the multi-core CPU. Then, we replaced the MUMPS solver with the CRVPINN solver - we used a single NVIDIA A100-SXM4-40GB GPU, pretrained for the initial state using 20,000 iterations (which takes 500 seconds). After that, we ran the training every 10 iterations of the IGA-ADS solver (we execute 100 iterations of the CRVPINN solver to compute the update of the pressure, which costs around 3-5 seconds). The total CRVPINN solver execution time on average is around $500 + 100 \times 4  = 900$ seconds. 
This is presented in Figure \ref{fig:IGACRVPINN}.
For example, for the mesh dimensions of $128^2$ for the IGA-ADS solver and the MUMPS solver, compared with $100^2$ collocation points for the CRVPINN solver, the hybrid IGA-ADS solver coupled with the MUMPS solver takes a total of around 1200 seconds (saturation integration) plus 3800 seconds (pressure integration) plus 5800 seconds (pressure solver), totaling 10,800 seconds (see Figure \ref{fig:IGAMUMPS}). The initialization and saturation solver times are negligibly small. On the other hand, the IGA-ADS solver coupled with the CRVPINN solver takes a total of around 1200 seconds (saturation integration) plus around 1300 seconds (communication between the IGA-ADS and CRVPINN solvers) plus 900 seconds (CRVPINN solver), totaling 3,400 seconds. We compare 10,800 seconds with 3,400 seconds. Our IGA-ADS solver, coupled with the CRVPINN solver, is over 3 times faster.

\section{Conclusions and Future Work}

In this work, we have demonstrated some accuracy and timing results of our implementation of the $CO_2$ sequestration hybrid solver using the  IGA-ADS and CRVPINN solver. 
The hybrid solver has shown capability in solving both simple (uniform permeability) and more complex (three realistic permeability maps) cases of the direct problem. 
It has also been able to handle the simulation well using a single node from the ARES supercomputer at ACK CYFRONET \cite{Wiatr}.
To be more specific, the $128^2$ mesh size simulation was performed on 12 cores and used 24 GB of RAM, while the CRVPINN was run on a single NVIDIA-A100 GPU. We have also tested the scalability of the IGA-ADS solver on extremely low resources (some test simulations were performed on a laptop with an AMD Ryzen 5 5500 U CPU @ 2.1 GHz and only 10 GB of RAM). The CRVPINN solver can also run on the free version of Google Colab, which allows for the usage of our solver even without access to state-of-the-art computational resources.
The sample results correspond to the expected behavior of the system, and the solver has been able to work with the parameter values given from the literature \cite{ebigbo_co2_2007,shokouhi_physics-informed_2021}.
We have compared our hybrid IGA-ADS solver, coupled with the CRVPINN method, with a baseline of the IGA-ADS solver coupled with the MUMPS direct solver. Our hybrid solver is over 3 times faster on the computational node from the ARES cluster of ACK CYFRONET.
%We have also conducted a comparison of the solver results with the parameters in SI and CGS units, and no visible differences were found. This serves as an additional validation of the solver.

%We plan extensive testing of the solver with various parameter configurations and prolonged simulation timescales in a short-term future.
%We are now finishing the setup of the test routines for running the code on the HPC center: ACK CYFRONET AGH.
Our future work will include the solution of the inverse problem: identifying the optimal location for the $CO_2$ injection well and the values of the adjustable parameters for a given reservoir. 
In the meantime, experiments on the $H_2$ storage problem may also be conducted.
As for the computational cost, we plan to work on lowering the cost of the communication between the IGA-ADS and CRVPINN solvers and further parallelizing the CRVPINN solver on several GPUs.
%After the completion of the direct finite element solver in the IGA-ADS framework, we plan to implement this problem using Physics-Informed Neural Networks (PINNs) \cite{karniadakis_physics-informed_2021}. 
%We will use recent Maczuga et al.'s work as a starting point for that step \cite{maczuga_physics_2024}.

Further improvement of the simulation quality can be obtained using $hp$-adaptive methods \cite{FI1,FI2}.
Future work may involve solving the inverse problem related to identifying the optimal location for pumps. They can be solved using multi-deme, twin adaptive strategy methods  \cite{Barabasz01012011,NaturalComputing}.

\section{Acknowledgements}

The work of Askold Vilkha and Maciej Paszyński has received funding from the European Union's Horizon Europe research and innovation programme under the Marie Sklodowska-Curie grant agreement No 101119556. We gratefully AcKnowledge the Polish high-performance computing infrastructure PLGrid (HPC Center: ACK Cyfronet AGH) for providing computer facilities and support within computational grants no. PLG/2024/017915, PLG/2024/017764, PLG/2025/018925.

\section*{Declaration of Generative AI and AI-assisted technologies in the writing process}
During the preparation of this work, the author(s) used Gemini (Google)  and ChatGPT to assist with writing, correcting English language usage, and improving text clarity and organization. After using this tool/service, the author(s) reviewed and edited the content as needed and take(s) full responsibility for the content of the publication.

\bibliographystyle{unsrt}
\bibliography{references}

\end{document}